\documentclass[12pt,reqno]{article}
mm \textheight=8.9in

\usepackage{amssymb}
\usepackage{amsmath}
\usepackage{amsthm}

\usepackage{color}
\usepackage{pb-diagram}
\newcommand{\tw}[3]{{$#1$}${\,\scriptscriptstyle {#2}}\atop\raise9pt\hbox{$\scriptstyle\tp$} ${$#3$}}
\newcommand{\st}[1]{\mbox{${\,\scriptscriptstyle {#1}}\atop\raise5.5pt\hbox{$*$}$}}
\newcommand{\braid}[3]{{#1}$\lower4pt\hbox{$\tp\atop\raise4pt
            \hbox{$\scriptscriptstyle {#2} $}$}${#3}}
\newcommand{\Us}[3]{{{}_{#1} \hspace{-0.1pt}\check \U_{#2}}^{#3}}
\newcommand{\btr}{\raise1.2pt\hbox{$\scriptstyle\blacktriangleright$}\hspace{2pt}}
\newcommand{\id}{\mathrm{id}}

\newcommand{\tp}{\otimes}

\newcommand{\Lc}{\mathcal{L}}
\newcommand{\A}{\mathcal{A}}

\newcommand{\Hg}{\mathcal{H}}
\newcommand{\Rg}{\mathcal{L}}
\newcommand{\Bg}{\mathcal{B}}

\newcommand{\Kg}{\mathcal{K}}
\newcommand{\Tg}{\mathcal{T}}
\newcommand{\B}{\mathcal{B}}
\newcommand{\U}{\mathcal{U}}
\newcommand{\D}{\mathfrak{D}}

\newcommand{\Fg}{\mathfrak{F}}

\newcommand{\J}{\mathcal{J}}
\newcommand{\N}{\mathbb{N}}
\newcommand{\Cc}{\mathcal{C}}

\renewcommand{\S}{\mathfrak{S}}
\newcommand{\Ha}{\mathcal{H}}

\newcommand{\Q}{\mathcal{Q}}

\newcommand{\F}{\mathcal{F}}
\newcommand{\Ru}{\mathcal{R}}
\newcommand{\ve}{\varepsilon}
\newcommand{\gm}{\gamma}
\newcommand{\Gm}{\Gamma}

\newcommand{\e}{\epsilon}

\newcommand{\la}{\lambda}

\newcommand{\End}{\mathrm{End}}

\newcommand{\Aut}{\mathrm{Aut}}

\newcommand{\Hom}{\mathrm{Hom}}

\newcommand{\Tr}{\mathrm{Tr}}
\newcommand{\tr}{\triangleright}

\newcommand{\btl}{\mbox{\raise1.1pt\hbox{$\scriptstyle\blacktriangleleft$}}}

\newcommand{\nn}{\nonumber}

\newcommand{\K}{\mathcal{K}}

\newcommand{\si}{\sigma}
\newcommand{\al}{\alpha}
\newcommand{\vt}{\vartheta}
\newcommand{\zt}{\zeta}
\newcommand{\bt}{\beta}
\newcommand{\Lcc}{{\scriptscriptstyle \mathcal{L}}}

\newcommand{\be}{\begin{eqnarray}}
\newcommand{\ee}{\end{eqnarray}}

\newtheorem{thm}{Theorem}[section]
\newtheorem{propn}[thm]{Proposition}
\newtheorem{lemma}[thm]{Lemma}

\theoremstyle{definition}
\newtheorem{remark}[thm]{Remark}
\newtheorem{definition}[thm]{Definition}
\newtheorem{example}[thm]{Example}

\newcount\prg

\newcommand{\parag}{\advance\prg by1 {\noindent\bf\thesection.\the\prg\hspace{6pt}}}
\newcommand{\select}[1]{\textcolor{red}{\bf\em #1}}
\begin{document}
\title{Dynamical reflection equation\footnote{
This research is supported in part
by the Israel Academy of Sciences grant no. 8007/99-03,
the Emmy Noether Research Institute for Mathematics,
the Minerva Foundation of Germany,  the Excellency Center "Group
Theoretic Methods in the study of Algebraic Varieties"  of the Israel
Science foundation, and by the RFBR grant no. 03-01-00593.}
}
\author{P. P. Kulish$^\dag$
\hspace{3pt} and A. I. Mudrov$^{\ddag,\natural}$
\\[15pt]
\em \small To the memory of Joseph Donin}
\date{}
\maketitle
\hspace{20pt}\vbox{\noindent \small
$^\dag$St.-Petersburg Department of Steklov Mathematical Institute,
Fontanka 27,
\\
\vbox{191011 St.-Petersburg, Russia, e-mail: \texttt{kulish@pdmi.ras.ru}}
\\
$^{\ddag}$Department of Mathematics, Bar Ilan University, 52900 Ramat Gan,
Israel,
\\
\vbox{e-mail: \texttt{mudrova@math.biu.ac.il}}
\\
$^{\natural}$Max-Planck Institut f$\ddot{\rm u}$r Mathematik, Vivatsgasse 7, D-53111 Bonn, Germany.
\\
\vbox{e-mail: \texttt{mudrov@mpim-bonn.mpg.de}}
}

\begin{abstract}
We construct a dynamical reflection equation algebra, $\tilde {\mathcal{K}}$, via a dynamical twist of the ordinary reflection
equation algebra. A dynamical version of the reflection equation
is deduced as a corollary.
We show that  $\tilde {\mathcal{K}}$ is a right comodule
algebra over a dynamical analog of the
Faddeev-Reshetikhin-Takhtajan algebra equipped with a structure of
right bialgebroid. We introduce dynamical trace and use it for
constructing central elements of $\tilde {\mathcal{K}}$.
\end{abstract} {\small \underline{Key words}:
Dynamical Yang-Baxter equation, dynamical reflection equation, quantum groupoids.\\
\underline{AMS classification codes}: 17B37, 81R50.}
\maketitle
\tableofcontents
\section{Introduction}
The dynamical Yang-Baxter equation (dYBE)
appeared in integrable models of conformal field theories
\cite{AF,ABB,Fad,F,GN}. It is a generalization of the ordinary
Yang-Baxter equation (YBE), which is the most important equation
of the quantum inverse scattering method. In the simplest version of
abelian base, the dYBE and related structures were studied
in \cite{EV1,EV2,ES2,S}, see also the lecture course \cite{ES2} and references therein.
The general formulation of dYBE (over arbitrary
base) was motivated by its relation to equivariant quantization
on G-spaces and developed  in \cite{DM1}; that formulation was based on
ideas of \cite{EV1,EV2,Xu2}.

The theory of the YBE may be considered as a part of the Hopf algebra
(more generally,  bialgebra) theory, \cite{Dr1}. The natural analog of Hopf algebra for
the dYBE is bialgebroid. The remarkable distinction of bialgebroids from Hopf algebras is
replacement of the field of scalars by an algebra that is non-commutative in general;
it is called "base algebra".
Roughly speaking, the connection between the theory of dYBE and bialgebroids
is established via the base algebra, whose spectrum is set to be
the space of dynamical parameter.

Many standard concepts of the Hopf algebra theory, such as tensor product of
modules, twists, module and comodule algebras can be formulated for bialgebroids,
\cite{Szl,Xu1,Sch}.
One can also define a quasitriangular structure on bialgebroids,
which in special cases is expressed through a dynamical R-matrix, \cite{DM2}.
In the present paper, we study dynamical analogs of two important
objects of the Hopf algebra theory: the Faddeev-Reshetikhin-Takhtajan (FRT)
and reflection equation (RE) algebras, \cite{FRT,KSkl,KS}.
The latter is given our special consideration.

A naive point of view on the dRE is to take for it the equation satisfied by
the squared permutation in the dynamical category, \cite{DM1}, as it
goes in the Hopf algebra theory. The question is what are the properties
of that equation and what algebraic structures are behind it.
The presented analysis of those structures supports that "naive" point of view
and shows that the theory of the dRE is analogous to
the theory of its non-dynamical counterpart.

Below we obtain the following results.
We develop a theory of dynamical reflection equation over arbitrary base
(by a base we understand a pair of a Hopf algebra and its base algebra, see \cite{DM1}).
In particular, we associates to any matrix solution $\tilde R$ the dYBE over a base algebra
$\Lc$ a
right $\Lc$-bialgebroid called dFRT algebra. We define
a dRE algebra as a right comodule algebra over the dFRT bialgebroid.
We concentrate our study on the situation when $\tilde R$ is
the image of the universal dynamical R-matrix $\tilde \Ru$ obtained
by a dynamical twist from a constant universal R-matrix, $\tilde \Ru=\F^{-1}_{21}\Ru\F$.
Under this assumption, we prove that the dFRT and dRE algebras
are modules over certain quantum groupoids and can be obtained by
dynamical twists from the FRT and RE algebras.
For twisted dynamical R-matrices, we introduce  dynamical trace
and prove that traces of powers of the dynamical RE matrices
lie in the center of the dRE algebra.

Let us mention two papers that are relevant to the subject in question.
A version of quadratic algebra depending on a dynamical parameter
corresponding to the abelian base was proposed in \cite{NAR}.
That algebra is different from what is studied in the present paper.
It is defined by four different matrices, whose precise
meaning has to be understood (except for one, which is a solution of the dYBE).
Another version of dynamical quadratic algebra was suggested in \cite{FHS}.
Essentially, that is a particular case of our dRE algebra specialized for
the abelian base.

The article is organized in the following way.
\begin{itemize}
\item
Section \ref{secHA} recalls the definition of left and right bialgebroids
 and gives some examples.
\item
Section \ref{secdFRTREa} introduces dFRT and dRE algebras for any dynamical permutation $\tilde S$.
\item
Section \ref{secTwFRTREa} discusses a way of constructing the dRE algebra by dynamical twist,
 following an analogy with the Hopf algebra case.
\item
Section \ref{secASDT} develops a technical machinery necessary for realization of the
 strategy adopted in Section \ref{secTwFRTREa}.
\item
Section \ref{secGTDT} studies symmetries of dynamical twists and dynamical R-matrices.
\item
Section \ref{secDTREa} constructs the dFRT and dRE algebras via dynamical twists
 of the ordinary FRT and RE algebras.
\item
Section \ref{secDTC} introduces dynamical trace used for constructing central
elements of the dRE algebra (in the twisted case).
\end{itemize}
\vspace{6pt}
{\bf Notation and conventions}\\
All the vector spaces are over the ground field $k$ of zero characteristic.
Algebras are associative unital algebras over the ground field. The subscript ${op}$ denotes
the opposite multiplication.
To make formulas more readable, we denote the inverse by bar, $\bar a:=a^{-1}$.
Coproduct, counit and antipode are denoted by $\Delta$, $\ve$, and $\gm$
without reference to a particular Hopf algebra or bialgebroid if they are clear from the context.
The symbol $\tp$ stands for the tensor product over the ground field. Tensors products over
other rings are indicated  explicitly.
We use the standard Sweedler notation for the coproduct $\Delta(x)=x^{(1)}\tp x^{(2)}$.
Similar convention is used for a Hopf algebra coaction, say, $\delta(v)=v^{(1)}\tp v^{[\infty]}$,
where the factor of the same nature as $v$ is marked with $[\infty]$ for left- and with $[0]$
for right comodules. The iterated coproducts are denoted by $\Delta^2:=\Delta$, $\Delta^3:=(\Delta\tp\id)\circ \Delta$ and so on.
Similar notation is used for iterated coactions.
Factors of tensor objects are labelled by positive integers in the standard way.
Notation $F_{(12)3}$ or $F_{1(234)}$  {\em etc }  means $(\Delta\tp \id)(F)$ or
$(\id\tp \Delta^3)(F)$ {\em etc } for an element $F$ from the tensor square of a Hopf algebra.
More special notation will be introduced in Subsection \ref{subsecSDN}.

\section{Bialgebroids}
\label{secHA}

The concept of bialgebroid generalizes that of bialgebra.
Bialgebroids naturally arise in the theory of dYBE as
bialgebras (Hopf algebras) in the theory of YBE.
Similarly to bialgebras,
representations of a bialgebroid form a monoidal category.
The distinction from the bialgebra case is that such monoidal
category has no fiber functor to the vector category. Rather,
it has a fiber functor to the monoidal category of $\Lc$-bimodules,
where $\Lc$ is an algebra called base.

Below we give, for further convenience, definitions
of left and right bialgebroids. For details, the reader is referred
to \cite{Lu} and \cite{Szl}. Bialgebroids in connection with
the dYBE were studied in \cite{Xu1} and \cite{DM2}.

Suppose we are given with a pair of associative algebras $(\Lc,\Bg)$,
a homomorphism $s\colon \Lc\to \Bg$ and anti-homomorphism $t\colon \Lc \to \Bg$
called source and target maps, respectively. Assume that their images commute
in $\Bg$.
One can consider the following two $\Lc$-bimodule structures on $\B$ induced by these maps.
One of them is induced by the left regular representation,
$
\la.b.\mu:=s(\la)t(\mu)b
$,
$\la,\mu\in \Lc$, $b\in \Bg$.
The other comes from the right regular action,
$
\la.b.\mu:=b t(\la)s(\mu)
$,
$\la,\mu\in \Lc$, $b\in \Bg$.
\subsection{Left bialgebroids}
\label{subsecLBA}
In this subsection we fix the $\Lc$-bimodule structure on  $\Bg$ defined via the left regular action.
The tensor square $\Bg\tp_\Lc \Bg$ is not an algebra, in general. However, it contains
a natural algebra as a subspace.
Namely,
put
\be
\Bg\times_\Lc \Bg=\{z\in \Bg\tp_\Lc \Bg\;| \>z\bigl( t (\la)\tp 1\bigr) = z\bigl(1\tp  s (\la)\bigr), \forall \la\in \Lc\}.
\label{subalg}
\ee
Clearly, $\Bg\times_\Lc \Bg$ is a unital associative algebra.
\begin{definition}
\label{Lu}
The quadruple $(\Bg,\Lc,s,t)$ is called
a \select{left bialgebroid over  base} $\Lc$ if there exist
\begin{enumerate}
\item a coassociative $\Lc$-bimodule map
(comultiplication) $\Delta\colon \Bg\to \Bg\times_\Lc \Bg$
that is a homomorphisms of associative algebras.
\item a bimodule map
(counit) $\ve\colon \Bg\to \Lc$ such that $\ve(1_\Bg)=1_\Lc$,
\be
&\ve\bigl(a\; (s\circ\ve) (b)\bigr)=\ve(ab)=\ve\bigl(a\; (t\circ\ve) (b)\bigr), \quad a,b\in \Bg,
\quad \mbox{ and}
\label{ve_l}
\\[7pt]
&(\ve\tp_\Lc\id_\Bg)\circ \Delta =\id_\Bg= (\id_\Bg\tp_\Lc\ve)\circ
\Delta \label{counit}
\ee
under the identification
$\Lc\tp_\Lc\Bg\simeq \Bg\simeq  \Bg\tp_\Lc\Lc$.
\end{enumerate}
\end{definition}

Condition 2 implies the identities $\ve\circ s=\ve\circ t=\id_\Lc$ and
makes $\Lc$ a left $\Bg$-module by
\be
a\vdash\la:=\ve\bigl(a s(\la)\bigr)=\ve\bigl(a t(\la)\bigr),
\quad
a\in \Bg, \la\in \Lc,
\label{anchor_l}
\ee
where the right equality is a consequence of (\ref{ve_l}).
The $\Lc$-bimodule structure on $\Lc$ induced by this action coincides with
the standard one. The action (\ref{anchor_l}) is called (left) \select{anchor}.
In fact, the anchor determines the counit together with condition
(\ref{ve_l}) by $\ve(a)=a \vdash 1_\Lc$ for $a\in \Bg$.

Any left $\Bg$-module $V$ is a natural $\Lc$-bimodule. This correspondence
is called the forgetful functor. Given two $\Bg$-modules $V$ and $W$,
the tensor product $V\tp_\Lc W$ acquires a left $\Bg$-module structure
via the coproduct, due to condition 1 of Definition \ref{Lu}.
The whole set of axioms from Definition \ref{Lu} ensures that the
left $\Bg$-modules form a monoidal category,
with $\Lc$ being the unit object. The forgetful functor  to
the category of $\Lc$-bimodules is strong monoidal, i.e. preserves tensor products.
Conversely, suppose there is a  monoidal structure on the category of left $\Bg$-modules.
Suppose the forgetful functor
to the category of $\Lc$-bimodules is strong monoidal.
Then $\Bg$ is a left $\Lc$-bialgebroid, see \cite{Szl}.
\subsection{Right bialgebroids}
The theory of right bialgebroids is fully parallel to the theory of left ones.
Let us write out the definition for future convenience.

Now assume the alternative $\Lc$-bimodule structure on  $\Bg$, that is specified by the right regular action.
The algebra $\Bg\times_\Lc \Bg$ is defined in this case as
\be
\Bg\times_\Lc \Bg=\{z\in \Bg\tp_\Lc \Bg\;| \>\bigl( s (\la)\tp 1\bigr)z = \bigl(1\tp  t (\la)\bigr)z, \forall \la\in \Lc\},
\label{subalg_r}
\ee
\begin{definition}
\label{Lu_r}
The quadruple $(\Bg,\Lc,s,t)$ is called
a \select{right bialgebroid over  base} $\Lc$ if there exist
\begin{enumerate}
\item a coassociative $\Lc$-bimodule map homomorphism
(comultiplication) $\Delta\colon \Bg\to \Bg\times_\Lc \Bg$
that is a homomorphisms of associative algebras.
\item a bimodule map
(counit) $\ve\colon \Bg\to \Lc$ such that $\ve(1_\Bg)=1_\Lc$,
\be
&\ve\bigl((s\circ\ve) (a)\; b\bigr)=\ve(ab)=\ve\bigl( (t\circ\ve) (a)\;b\bigr), \quad a,b\in \Bg,
\quad \mbox{ and}
\label{ve_r}
\\[7pt]
&(\ve\tp_\Lc\id_\Bg)\circ \Delta =\id_\Bg= (\id_\Bg\tp_\Lc\ve)\circ
\Delta \label{counit_r}
\ee
under the identification
$\Lc\tp_\Lc\Bg\simeq \Bg\simeq  \Bg\tp_\Lc\Lc$.
\end{enumerate}
\end{definition}
Similarly to the left bialgebroids, one can define the (right) anchor action setting
\be
\la\dashv a:=\ve\bigl(s(\la)a\bigr)=\ve\bigl( t(\la)a\bigr),
\quad
a\in \Bg, \la\in \Lc.
\label{anchor_r}
\ee
The anchor determines the counit by $\ve(a)=1\dashv a$ for $a\in \Bg$.

The right modules over right bialgebroids form a monoidal category with the fiber
functor to the monoidal category of $\Lc$-bimodules.
\begin{remark}
It is possible to define twist and quasitriangular structure for bialgebroids,
as well as for Hopf algebras. The reader is referred to \cite{Xu1} and \cite{DM2} for details.
\end{remark}

\subsection{Some examples}
Let us give some examples of bialgebroids that arise in the theory of dYBE.
More examples will appear bellow.
Bialgebroids related to the dYBE are defined over a base that is commutative
with respect to some permutation.
More specifically, let us fix  a Hopf algebra $\Hg$ with the coproduct $\Delta$,
and invertible antipode $\gm$.  Recall from \cite{DM1} that
an $\Hg$-base algebra $\Lc$ is a left $\Hg$-module algebra and left $\Hg$-comodule algebra; it is
a Yetter-Drinfeld module with respect to these structures and "commutative" in the following sense:
\be
\label{q-com}
\la\mu = (\la^{(1)}\tr\mu) \la^{[\infty]}, \quad \la, \mu \in \Lc.
\ee
Here $\la\mapsto \la^{(1)}\tp \la^{[\infty]}\in \Hg\tp \Lc$ denotes the coaction $\delta$ and $\tr$ stands
for the $\Hg$-action on $\Lc$.
Algebraically, $\Lc$ can be defined as a module algebra over the double $\D\Hg=\Hg\bowtie\Hg^*_{op}$,
\cite{Dr1},
that is commutative with respect to the standard quasitriangular structure $\Theta\in \Hg^*_{op}\tp \Hg\subset(\D\Hg)^{\tp 2}$.
The coaction may be written in terms of $\Theta$ and the $\D\Hg$-action
(for this we reserve notation $\tr$ throughout the whole paper)
\be
\label{Theta_delta}
\delta\colon\la\mapsto \Theta_2\tp \Theta_1\tr\la.
\ee
Note that being a Yetter-Drinfeld module over $\Hg$ is immediate from this
representation and the intertwining axiom for R-matrices.
It is convenient to use representation (\ref{Theta_delta}), although all the results would be true without
this simplifying assumption.
\begin{example}
\label{ex_bialg1}
Consider the smash product $\Lc\rtimes \D\Hg$ algebra with the multiplication
\be
(\la\tp f)(\mu \tp g)&:=&\la (f^{(1)}\tr \mu) \tp f^{(2)}g,
\quad
\la,\mu\in \Lc,
\quad f,g\in \Hg.
\ee
The algebra  $\Lc\rtimes \D\Hg$ is a left $\Lc$-bialgebroid
with the subalgebra $\Lc\rtimes \Hg$ being a sub-bialgebroid.
The construction is due to \cite{Lu} and goes as follows.
The source and target maps are set to be
$$
s(\la)=\la\tp 1,
\quad
t(\la)=\bar \Theta_1\tr\la \tp \bar \Theta_2
$$
The counit and coaction on $\Lc\rtimes \D\Hg$ are defined as
\be
\ve(\la\tp h)&=&\la\ve_\Hg(h),
\\
\Delta(\la\otimes h)&=& (\la\otimes  h^{(1)})\tp_\Rg(1\otimes h^{(2)}).
\ee
The anchor action is given by $(\la\tp h)\vdash \mu=\la(h\tr\mu)$, for $\la\tp h\in \Lc\rtimes\D\Hg$ and $\mu\in \Lc$.
\end{example}
\begin{example}
\label{ex_bialg2}
Let $\Lc$ be an $\Hg$-base algebra with the $\Hg$-action $\tr$ and the coaction $\delta$.
Then $\Lc_{op}$ is a base algebra over $\Hg_{op}$ with respect to
the $\Hg_{op}$-action $x\tp\ell\mapsto \gm^{-1}(x)\tr \ell$ and
the same coaction $\delta$ considered now as a map $\Lc_{op}\to \Hg_{op}\tp \Lc_{op}$.
Along the line of Example \ref{ex_bialg1}, one can construct the left $\Lc_{op}$-bialgebroids
$\Lc_{op}\rtimes \D\Hg_{op}$ and $\Lc_{op}\rtimes \Hg_{op}$.
\end{example}

Similarly to Hopf algebras, one can define biideals in bialgebroids.
Namely, an $\Lc$-bimodule $\J\subset \Bg$
is called biideal if $\Delta(\J)\subset \J\tp_\Lc \Bg+\Bg\tp_\Lc \J$
and $\ve(\J)=0$. The latter condition holds if and only if $\J$ lies in the kernel
of the anchor. The quotient of a bialgebroid by a biideal is a bialgebroid.
Conversely, the kernel of a bialgebroid homomorphism is a biideal.

\begin{example}
\label{ex_bialg3}
The bialgebroid $\Lc\rtimes \D\Hg$ from \ref{ex_bialg1} is not quasitriangular.
Here is  an example of quasitriangular bialgebroid (quantum groupoid) from \cite{DM2}.
The quotient $\D\Hg_\Lc$ of  $\Lc\rtimes \D\Hg$ by
the two-sided ideal generated
by $\bar \Theta_1\tr\la \tp \bar \Theta_2- \Theta_2\tr\la \tp \Theta_1$, $\forall \la\in \Lc$,
is again a left $\Lc$-bialgebroid (recall that bar means the inverse). It is quasitriangular,
with the universal R-matrix obtained by a natural projection from the universal R-matrix of
$\D\Ha$, see \cite{DM2} for details.
\end{example}

\section{Dynamical FRT and RE algebras}
\label{secdFRTREa}

Let $V$ be a finite dimensional $\Hg$-module and $\rho$ denote the homomorphism
from $\Hg$ to $\End(V)$.
Let
$\{e^i_j\} \subset \End(V)$ denote the standard matrix base with the multiplication
to $e^i_j e^l_k = \delta^l_j e^i_k$; here $\delta^l_j$ is the
Kronecker symbol.

Let $\tilde S$ be an element from $\End^{\tp2}(V)\tp \Lc$ that is invariant in the following sense:
\be
\label{H-inv}
h\tr \tilde S = (\rho^{\tp 2}\circ\Delta)(\gm h^{(1)})\tilde S (\rho^{\tp 2}\circ\Delta)(h^{(2)}) \quad \mbox{for all }h\in \Hg.
\ee
It follows from this transformation law that $\tilde S$ commutes with
$
(\rho^{\tp 2}\circ\delta^2)(\la)
$
for all $\la\in \Lc$ (here $\delta^2$ is the two-folded coaction).
Put $\tilde S_{op}:= \tilde S_1\rho(\gm \Theta^{(2)}_2) \tp \tilde S_2\rho(\gm\Theta^{(1)}_2)\tp \Theta_1\tr \tilde S_\Lcc$.

In the next two subsections we associate with the matrix $\tilde S$ the dynamical analogs of two important algebras
in the theory of quantum groups, the FRT and RE algebras.

\subsection{dFRT algebra associated with a dynamical matrix $\tilde S$}
We call an algebra  $\A$ in the monoidal category of $\Lc$-bimodules
an $\Lc$-bimodule algebra.
It is an associative algebra whose multiplication and $\Lc$-bimodule
structure are compatible, namely
$
(a.\la)b= a (\la.b)
$
for all $a,b\in \A$ and $\la\in \Lc$. A natural example of an $\Lc$-bimodule algebra is obtained
when the bimodule structure is induced by a homomorphism $\Lc\to \A$. This is the case
if and only if the left and right $\Lc$-actions satisfy $\la.1_\A=1_\A.\la$ for all $\la\in \Lc$.

Let $\{T^i_j\}\subset \End^*(V)$ be the dual base to $\{e_i^j\}$.
Consider the space  $\End^*(V)$ as a left $\Hg_{op}\tp \Hg$-module with respect to
the action
\be
 x\btr T:= T \rho(x), \quad
y\btr T:= \rho(y) T,
\label{bi-act}
\ee
where $x\in \Hg$,
$y\in \Hg_{op}$, and $T:=\sum_{i,j} e_i^j\tp T^i_j$.
The action of $\Hg_{op}\tp \Hg$ naturally extends to an action on the free algebra
$k\langle T^i_j\rangle:=k\langle \{T^i_j\}\rangle$
making the latter an $\Hg$-bimodule algebra.
Consider the smash product $\Fg(T):=k\langle T^i_j\rangle \rtimes (\Lc_{op}\tp \Lc)$
determined by the permutation relations
\be
\label{FRTbimod}
\la T = T\rho(\la^{(1)})
\>\la^{[\infty]} , \quad  \mu T   = \rho(\mu^{(1)})T\> \mu^{[\infty]}
\ee
for $\la \in \Lc$ and $\mu \in \Lc_{op}$. It is easy to see that $\Fg(T)$
is a left module algebra over the Hopf algebra
$\bigl(\mbox{\tw{\Hg_{op}}{\bar\Theta}{\D\Hg_{op}}}\bigr)\tp \bigl(\mbox{\tw{\Hg}{\Theta}{\D\Hg}}\bigr)$,
where $\D\Ha$ and $\D\Ha_{op}$ act essentially on the $\Lc$ and $\Lc_{op}$ factors, respectively.
Here \tw{\:}{\F}{\:} denotes the twisted tensor product of two Hopf algebras by a bicharacter $\F$,
see \cite{RS} and also \cite{DM3}.
We  also regard $\Fg(T)$ as an $\Ha$-bimodule by the
Hopf algebra embedding
$$
\Ha_{op}\tp \Ha\stackrel{\Delta\tp \Delta}{\longrightarrow}
\bigl(\mbox{\tw{\Hg_{op}}{\bar\Theta}{\D\Hg_{op}}}\bigr)\tp \bigl(\mbox{\tw{\Hg}{\Theta}{\D\Hg}}\bigr).
$$

Consider the ideal $\J_{dFRT}$ in $\Fg(T)$ generated by the relations
\be
\label{dFRTrel}
 \tilde S_{op} T_1T_2 = T_1T_2 \tilde S,
\ee
where the matrix coefficients of  $\tilde S_{op}$ and $\tilde S$ belong respectively to $\Lc_{op}$ and $\Lc$.
\begin{definition}
\label{dFRTa}
The dFRT algebra $\tilde \Tg$ associated with the matrix $\tilde S$ is
the quotient of $\Fg(T)$ by the ideal
$\J_{dFRT}$.
\end{definition}

It is straightforward to check, using (\ref{H-inv}), that
$\J_{dFRT}$ (and therefore $\tilde \Tg$) is an $\Hg$-bimodule.
The natural embedding
$$
(t\tp s)\colon\Lc\tp \Lc \to k\langle T^i_j\rangle\rtimes \bigl(\Lc_{op}\tp \Lc\bigr)
$$
makes $\Fg(T)$, and therefore $\tilde T$, an $\Lc$-bimodule via the right regular action.
These $\Ha$- and $\Lc$-bimodule structures on $\Fg(T)$ and $\tilde \Tg$
amount to a structure of  left module over the $\Lc_{op}\tp \Lc$-bialgebroid $(\Lc_{op}\rtimes \Hg_{op})\tp (\Lc\rtimes \Hg)$;
the latter is a tensor product of $\Lc_{op}$- and $\Lc$-bialgebroids,
see Examples \ref{ex_bialg1} and \ref{ex_bialg2}.

Our next objective is to introduce a structure of right bialgebroid on $\tilde \Tg$.
Let us start with the anchor action and the counit.
\begin{lemma}
The formulas
\be
\la\dashv \mu:=\la\mu
, \quad
\la\dashv \nu:=\nu\la
, \quad
\la\dashv T^i_j:=\rho(\la^{(1)})^i_j\la^{[\infty]},
\quad \la,\mu\in \Lc,
\quad \nu\in \Lc_{op}.
\label{anchor}
\ee
 define a right action of $\tilde \Tg$ on $\Lc$.
 \label{lemma_dFRTanch}
\end{lemma}
\begin{proof}
Formulas (\ref{anchor}) define an action of the free algebra $k\langle T^i_j\rangle$
and actions of $\Lc$ and $\Lc_{op}$. Thus, to
prove the statement, it suffices to check that this action respects
permutation rules of $\Lc$ and $\Lc_{op}$ with the generators $T^i_j$ and
annihilates the ideal $\J_{FRT}$. This is done by a straightforward calculation and left
to the reader as an exercise. To give a hint, in the  matrix form the action of $T$ can be written as
$\la\dashv T=\delta(\la)$ (we suppressed $\rho$). Then the relations (\ref{dFRTrel}) will reduce
to the equation
$\tilde S\delta^2(\la)-\delta^2(\la)\tilde S =0$, which holds by the $\Ha$-invariance assumption (\ref{H-inv}).
\end{proof}
The action (\ref{anchor}) induces an $\Lc$-bimodule structure on $\Lc$
by
$\al.\la.\bt= \la\dashv s(\bt)t(\al)$. By construction, this $\Lc$-bimodule structure coincides with the standard
one.

\begin{lemma}
The correspondence
\be
T^i_j\mapsto\sum_{k}T^k_j\tp_\Lc T_k^i
 \quad\mbox{ or }\quad
T\mapsto T^{(1)}T^{(2)}
\mbox{ in the matrix form}
\label{dFRTcopr}
\ee
defines a homomorphism of algebras $\Delta\colon \Fg(T) \to \Fg(T)\times_\Lc\Fg(T)$,
which is a homomorphism of $\Lc$-bimodules.
\label{lemma_dFRTcopr}
\end{lemma}
\begin{proof}
First of all notice that $\bigl(s(\la)T^{(1)}\bigr)T^{(2)}=T^{(1)}\bigl(t(\la)T^{(2)}\bigr)$,
as follows from (\ref{FRTbimod}). Therefore $\Delta(T^i_j)$ belong to
$\Fg(T)\times_{\Lc} \Fg(T)$. Extend $\Delta $ as a homomorphism to the free subalgebra $k\langle T^i_j\rangle \subset\Fg(T)$.
Extend $\Delta$  to $\Lc\subset \Fg(T)$ and $\Lc_{op}\subset \Fg(T)$
as an $\Lc$-bimodule map:
$$\Delta(\la)=1_{\Fg(T)}\tp_\Lc s (\la)
\quad\mbox{and} \quad
\Delta(\mu)=t (\mu)\tp_\Lc 1_{\Fg(T)}
$$
for
$\la\in \Lc$ and $\mu\in \Lc_{op}$. It is straightforward to check that the relations (\ref{subalg_r}) are preserved by $\Delta$.
Obviously $\Delta\bigl(\Fg(T)\bigr)\subset \Fg(T)\times_\Lc\Fg(T)$, so the lemma is proven.
\end{proof}

\begin{propn}
(i) The coproduct $\Delta$ from Lemma \ref{lemma_dFRTcopr} together with the anchor action (\ref{anchor})
define a right $\Lc$-bialgebroid structure on the $\Fg(T)$.
(ii) The ideal $\J_{dFRT}\subset\Fg(T)$ is a biideal, hence the projection from $\Fg(T)$ along $\J_{dFRT}$ makes
$\tilde \Tg$ a right $\Lc$-bialgebroid.
\end{propn}
\begin{proof}
Statement (i) is obvious. Indeed, $\Delta$ is apparently coassociative.
It is easy to see that the counit $\ve$ defined through the anchor by
\be
\ve(x):=1_\Lc\dashv x,\quad x\in \Fg(T),
\ee
satisfies (\ref{counit_r}).

Lemma \ref{lemma_dFRTanch} says that the ideal $\J_{dFRT}$ lies in the kernel of the
anchor action (\ref{anchor}). To prove (ii), we need to check consistency of the coproduct
(\ref{dFRTcopr}) with the multiplication in $\tilde \Tg$.
An easy verification that relations (\ref{FRTbimod})
are respected is left to the reader. We will check only the dFRT relations (\ref{dFRTrel}).
We have for $\Delta(\tilde S_{op}T_1T_2)$, modulo $\J_{dFRT}\tp_\Lc \Fg(T)+\Fg(T)\tp_\Lc \J_{dFRT}$:
\be
\bigl(t(\tilde S_{op})T^{(1)}_1T^{(1)}_2\bigr)(T^{(2)}_1T^{(2)}_2)
&=&
\bigl(T^{(1)}_1T^{(1)}_2s(\tilde S)\bigr) (T^{(2)}_1T^{(2)}_2)
=
(T^{(1)}_1T^{(1)}_2)\bigl(t(\tilde S_{op}) T^{(2)}_1T^{(2)}_2\bigr)
\nn
\\
&=&
(T^{(1)}_1T^{(1)}_2)\bigl(T^{(2)}_1T^{(2)}_2 s(\tilde S)\bigr).
\nn
\ee
Here we set $t(\tilde S_{op}):=(\id\tp t)(S_{op})$ and $s(\tilde S):=(\id\tp s)(S)$.
We employed relations (\ref{dFRTrel}), then we pulled the coefficients of the matrix $s(\tilde S)$ over to
the right $\tp_\Lc$-tensor factor to the rightmost position as the coefficients of the matrix $t(\tilde S)$ .
After that we pushed the entries of $t(\tilde S)$ to the left. This operation gave the matrix
factor $t(\tilde S_{op})$. Applying relations (\ref{dFRTrel}) once again we obtained the final expression, which is equal
to $\Delta(T_1T_2\tilde S)$.

Thus we proved that $\J_{dFRT}$ is a biideal
in $\Fg(T)$, hence the quotient $\tilde \Tg=\Fg(T)/\J_{dFRT}$ inherits an $\Lc$-bialgebroid structure.
\end{proof}
\begin{remark}
In the constructions of the present subsection, we never used the dYBE on $\tilde S$, but only the
$\Ha$-invariance (\ref{H-inv}). When $\Lc$ is the field of scalars, this construction
gives the bialgebra associated with arbitrary matrix $\tilde S$, not necessarily
as solution to the YBE, see \cite{FRT}.
\end{remark}

In what follows, we will need an extended dynamical FRT algebra whose description is given below.
It has the well known analog in the theory of quantum groups, \cite{FRT}.

\subsection{Extended dFRT algebra}
\label{subsecEdFRTa}
In this subsection we construct a right bialgebroid $\braid{\tilde \Tg}{\tilde S}{\tilde \Tg}$ generated
by the double set of generators $\{T^i_j\}$, $\{\bar T^i_j\}$ and its natural quotient bialgebroid.
We restrict ourselves  only with the description of $\braid{\tilde \Tg}{\tilde S}{\tilde \Tg}$
leaving verification of the bialgebroid axioms to the reader.

Instead of $k\langle T^i_j\rangle$, consider the free algebra $k\langle T^i_j,\bar T^i_j\rangle$
generated by the matrix elements of $||T^i_j||$ and $||\bar T^i_j||$  with
the following $\Hg\tp \Hg_{op}$-module structure:
\be
&x\btr T = T\rho(x)
,\quad
y\btr T   = \rho(y)T
,\quad
x \btr\bar T =\rho\bigl(\gm (x)\bigr) \bar T
, \quad
y\btr\bar T   =  \bar T \rho\bigl(\gm(y)\bigr),
\label{FRTanchor}
\ee
where $x\in \Hg$ and
$y\in \Hg_{op}$.
The rest of the construction replicates the construction of the algebra $\tilde\Tg$ in the previous
subsection.
Using the $\Hg_{op}\tp \Hg$-coaction on $ \Lc_{op}\tp\Lc$
and the $\Hg_{op}\tp \Hg$-module structure (\ref{FRTanchor}), we define
the smash product $\Fg(T,\bar T)=k\langle T^i_j,\bar T^i_j\rangle\rtimes (\Lc_{op}\tp\Lc)$.
There are the following commutation relations held in $\Fg(T,\bar T)$:
\be
\begin{array}{lllllll}
\la T &=& T\rho(\la^{(1)})\>\la^{[\infty]}
,
&
\la \bar T &=&\rho\bigl(\gm (\la^{(1)})\bigr) \bar T\>\la^{[\infty]}
,
\\
\mu T   &=& \rho(\mu^{(1)})T \>\mu^{[\infty]}
,\quad
&\mu\bar T &=& \bar T \rho\bigl(\gm(\mu^{(1)})\bigr) \>\mu^{[\infty]},
\end{array}
\label{TbarTL}
\ee
where $\la\in \Lc$ and $\mu\in \Lc_{op}$.
There exist obvious algebra and anti-algebra embeddings, $s$ and $t$, from
$\Lc$ to $\Fg(T,\bar T)$.
The algebra $\Fg(T,\bar T)$ is endowed with an $\Lc$-bimodule structure
via these maps and the right regular representation.
We consider the two-sided ideal in $\Fg(T,\bar T)$
generated by the relations
\be
&\tilde S_{op}T_1T_2=T_1T_2 \tilde S ,
\quad
T_2\tilde S\bar T_2 = \bar T_1 \tilde S_{op} T_1
,
\quad
\tilde S\bar T_2 \bar T_1=\bar T_2 \bar T_1 \tilde S_{op},
\label{T-barT}
\ee
and denote by $\braid{\tilde \Tg}{\tilde S}{\tilde \Tg}$ the quotient
of $\Fg(T,\bar T)$ by this ideal.
We define the right anchor action of $\Fg(T,\bar T)$ on $\Lc$ setting
\be
\begin{array}{rclrclc}
\la\dashv T^i_j & :=&\rho(\la^{(1)})^i_j\la^{[\infty]}, &\quad
\la\dashv \bar T^i_j&:=&\rho\bigl(\gm(\la^{(1)})\bigr)^i_j\la^{[\infty]},
\\
\la\dashv \mu &:=&\la\mu ,& \la\dashv \nu&:=&\nu\la
\end{array}
\label{anchor1}
\ee
for  $\la,\mu\in \Lc$ and $\nu\in \Lc_{op}$. We can prove that this action descends to
an action of the quotient algebra $\braid{\tilde \Tg}{\tilde S}{\tilde \Tg}$.

Then we introduce an $\Lc$-bimodule map from $\Fg(T,\bar T)$  to $\Fg(T,\bar T)\times_\Lc\Fg(T,\bar T)$ setting
it on the generators $\{T^i_j\}$ and $\{\bar T^i_j\}$ by
\be
\Delta(T^i_j)=\sum_{k}T^k_j\tp_\Lc T_k^i,
\quad
\Delta(\bar T^i_j)=\sum_{k}\bar T_k^i\tp_\Lc \bar T^k_j.
\ee
or, in the matrix form, by
\be
\Delta(T)=T^{(1)}T^{(2)},
\quad
\Delta(\bar T)=\bar T^{(2)}\bar T^{(1)}.
\ee
This correspondence  is extended to the entire $\Fg(T,\bar T)$ as an algebra homomorphism.
It is coassociative and defines, together with the anchor action (\ref{anchor1}),
a right $\Lc$-bialgebroid structure on $\Fg(T,\bar T)$.
The ideal of relations (\ref{T-barT}) is a biideal, hence the bialgebroid structure
is carried over to $\braid{\tilde \Tg}{\tilde S}{\tilde \Tg}$.

The algebra $\braid{\tilde \Tg}{\tilde S}{\tilde \Tg}$ admits a natural quotient by the
$\Ha$-invariant two-sided ideal specified by the relations
\be
\label{T-barT=1}
\bar T T = 1= T \bar T
\ee
in the concise matrix form. We call this quotient extended dFRT algebra and denote
by $\tilde\Tg_{ext}$.
In fact, relations (\ref{T-barT=1}) define a biideal, so $\tilde \Tg_{ext}$
is a right $\Lc$-bialgebroid.

\begin{remark}
When one puts $\Ha=\Lc=k$, the dFRT and dRE algebras $\tilde \Tg$ and $\tilde \Kg$ degenerate to the ordinary
FRT and RE algebras $\Tg$ and $\Kg$ associated with a matrix  $S\in \End(V\tp V)$.
\end{remark}
\subsection{Category of comodules over a bialgebroid}
As was mentioned in Subsection \ref{subsecLBA}, the category of modules over a bialgebroid
is equipped with a monoidal structure. Let us render the definition of the monoidal category of
comodules over a (right) bialgebroid following \cite{Sch}.
\begin{definition}
\label{comodule}
Let $\Bg$ be a right bialgebroid. An $\Lc$-bimodule $V$ is called right $\Bg$-comodule
if it is equipped with an $\Lc$-bimodule map (coaction) $\delta\colon V\to V\tp_\Lc \Bg$ fulfilling
\begin{enumerate}
\item
for all $\la\in \Lc$ and $v\in V$
\be
\la.v^{[0]}\tp_\Lc v^{(1)}=v^{[0]}\tp_\Lc t_\Bg(\la)v^{(1)}
,\quad \mbox{for}\quad
v^{[0]}\tp_\Lc v^{(1)}=\delta(v),
\label{comod_c1}
\ee
\item
coassociativity condition
\be
(\delta\tp_\Lc \id)\circ \delta=(\id \tp_\Lc \Delta)\circ \delta,
\ee
\item counital condition:
$(\id\tp_\Lc \ve_\Bg)\circ\delta\simeq \id$ under the identification $V\tp_\Lc\Lc\simeq V$.
\end{enumerate}
\end{definition}
All $\Bg$-comodules form a category, $\Cc$. The set $\Hom_\Cc(V,W)$ consists of maps
$V\stackrel{\phi}{\rightarrow}W$ such that the diagram
\be
   \dgARROWLENGTH=0.63\dgARROWLENGTH
\begin{diagram}
\node{V}\arrow{e,t}{\delta_V}\arrow{s,l}{\phi}
\node{V\tp_\Lc \Bg}\arrow{s,r}{\phi\tp_\Lc\id}
\\
\node{W}\arrow{e,t}{\delta_W}\node{W\tp_\Lc \Bg}
\end{diagram}
\nn
\ee
is commutative.

Let us introduce on $\Cc$ a structure of strict monoidal category.
Given two $\Bg$-comodules
$V,W$ define the $\Bg$-coaction on $V\tp_\Lc W$ as
$$v\tp_\Lc w\mapsto v^{[0]}\tp_\Lc w^{[0]} \tp_\Lc v^{(1)}w^{(1)}.$$
This map is correctly defined by virtue of (\ref{comod_c1})
and fulfills the axioms of coaction. The tensor product of two morphisms,
$\phi$ and $\psi$, is set to be $\phi\tp_\Lc \psi$. The unit object
in $\Cc$ is $\Lc$, with the coaction $\la\mapsto \la\tp_\Lc 1_\Bg$.
By construction, the category $\Cc$ has a forgetful functor to
the category of $\Lc$-bimodules, and this functor is strong monoidal.

Let $\A$ be an $\Lc$-bimodule algebra. Then
$$
\A\times_\Lc \Bg=\{z\in \A\tp_\Lc \Bg\>|\bigl(\la.\tp 1\bigr)z =\bigl(1\tp t_\Bg(\la)\bigr)z, \forall \la\in \Lc\},
$$
is again an $\Lc$-bimodule algebra. If $\A$ is a $\Bg$-comodule, then the coaction takes values in $\A\times_\Lc \Bg$.
An $\Lc$-bimodule algebra and $\Bg$-comodule is an algebra in the category $\Cc$ if and only if
the coaction $\delta\colon \A\to \A\times_\Lc \Bg$ is a homomorphism of $\Lc$-bimodule algebras.
\begin{definition}
An algebra in $\Cc$ is called $\Bg$-comodule algebra.
\end{definition}
Clearly the bialgebroid $\Bg$ itself is a $\Bg$-comodule algebra.

\subsection{dRE algebra associated with a dynamical matrix $\tilde S$}
\label{ssDFRTa}
This time denote the dual base to $\{e_i^j\subset \End(V)\}$ by $\{K^i_j\}\subset \End^*(V)$.
The left $\Hg$-action
\be
h\btr K= \rho\bigl(\gm(h^{(1)})\bigr)K \rho(h^{(2)}),
\quad
\mbox{where}
\quad
K:=\sum_{i,j} e_i^j\tp K^i_j,
\label{coadj}
\ee
on $\End^*(V)$  naturally extends to a left action on $k\langle K^i_j \rangle$
making it an $\Hg$-module algebra.
Consider the twisted tensor product of Hopf algebras \tw{\Hg}{\Theta}{\D\Hg} and
construct the corresponding twisted module algebra $\Fg(K):=k\langle K^i_j \rangle\rtimes \Lc$.
The algebra $\Fg(K)$ contains $\Lc$ as a subalgebra and thus it is a natural $\Lc$-bimodule.
The elements of $\Lc\subset \Fg(K)$ obey the following permutation rules with the generators
$K^i_j$:
\be
\la K = \rho\bigl(\gm(\la^{(1)})\bigr)K \rho(\la^{(2)})\> \la^{[\infty]}.
\label{Ktrans}
\ee
The actions of $\Hg$ and $\Lc$ on $\Fg(K)$ give rise to an action of the left $\Lc$-bialgebroid
$\Lc\rtimes \Hg$. In fact, $\Fg(K)$ is an algebra over $\Lc\rtimes \Hg$.
\begin{definition}
The dRE algebra associated with an $\Ha$-invariant element $\tilde S\in \End^{\tp 2}(V)\tp \Lc$ and
denoted further by  $\tilde \Kg$ is
the quotient of $\Fg(K)$ by the relations
\be
\label{dRE}
\tilde S K_2 \tilde S K_2=K_2 \tilde S K_2 \tilde S.
\ee
\label{defdREa}
\end{definition}
\vspace{-30pt}
Note with care that the coefficients of $\tilde S$ belong to $\Lc$ and do not commute with $K^i_j$.

One can check, using (\ref{H-inv}), that the
ideal in $\Fg(K)$ generated by (\ref{dRE}) is invariant with respect to $\Hg$,
so $\tilde \Kg$ is an $\Ha$-module algebra. It is also
is an $\Lc$-bimodule algebra whose $\Lc$-bimodule structure is induced by the homomorphism
$\Lc\to \tilde \Kg$ denoted further $s_{\tilde \Kg}$.
These two structures make $\tilde K$ a module algebra over the bialgebroid $\Lc\rtimes \Ha$.

\begin{propn}
The correspondence
$\delta_{\tilde \K}\colon K^i_j\tp \la\mapsto
\sum_{\al,\bt}(K^\al_\bt\tp 1_\Lc)\tp_\Lc (\bar T_j^\bt T^i_\al\tp 1_{\Lc_{op}}\tp \la)$
defines on $\tilde \Kg$ a structure of right $\braid{\tilde \Tg}{\tilde S}{\tilde \Tg}$-comodule algebra.
It also makes $\tilde \Kg$ a right $\tilde \Tg_{ext}$-comodule algebra via the projection
$\braid{\tilde \Tg}{\tilde S}{\tilde \Tg}\to \tilde \Tg_{ext}$.
\end{propn}
\begin{proof}
Put $\A:=\tilde \Kg$ and  $\B:=\braid{\tilde \Tg}{\tilde S}{\tilde \Tg}$ .
Let us show that the image of the "similarity transformation" $\delta_{\tilde \K}\colon K\mapsto \bar T K T $
lies in $\End(V)\tp \A\times_\Lc \Bg$.
Indeed,
\be
\bar T s_\A(\la) K T &=& \bar T \rho\bigl(\gm(\la^{(1)})\bigr) K \rho(\la^{(2)}) s_\A(\la^{[\infty]}) T
\nn\\
&=&
\bar T \rho\bigl((\gm(\la^{(1)})\bigr) K \rho(\la^{(2)})  T t_\Bg(\la^{[\infty]})
=
t_\Bg(\la)\bar T  K  T
.
\nn
\ee
The first and the last equalities are obtained using (\ref{Ktrans}) and (\ref{TbarTL}).
In the middle one, we employed the definition of tensor product over $\Lc$.

Next we prove that $\delta_{\tilde \K}$ defines an algebra homomorphism from $\A$ to $\A\times_\Lc\Bg$.
First of all, observe that
$s_\Bg(\la) (\bar T K T)= \rho(\gm \la^{(1)})(\bar T K T)\rho(\la^{(2)})\> s_\Bg(\la^{[\infty]})$.
Thus the relation (\ref{Ktrans}) is preserved.

Let us compute the expression
$s_\Bg(\tilde S) (\bar T_2 K_2 T_2 ) s_\Bg(\tilde S) (\bar T_2 K_2 T_2)$.
Using (\ref{dFRTrel}) twice, we find it equal to
\be
s_\Bg(\tilde S) (\bar T_2 K_2 T_2 ) s_\Bg(\tilde S) (\bar T_2 K_2 T_2)
&=&
s_\Bg(\tilde S) \bar T_2 K_2 \bar T_1 t_\Bg(\tilde S_{op}) T_1 K_2 T_2=
\nn\\
=s_\Bg(\tilde S) \bar T_2 \bar T_1  K_2 t_\Bg(\tilde S_{op}) K_2 T_1 T_2
&=&
\bar T_2 \bar T_1 t_\Bg(\tilde S_{op}) K_2 t_\Bg(\tilde S_{op}) K_2 T_1 T_2
.
\label{eq_aux1}
\ee
We can replace $t_\Bg(\tilde S_{op})$  by $s_\A(\tilde S)$ in the last expression. Indeed,
pull the entries of the right matrix $t_\Bg(\tilde S_{op})$ to
the right-most position. We have to commute them with the entries
of the matrices $T_1$ and $T_2$, according to (\ref{TbarTL}). Then we replace the entries of $t_\Bg(S)$ by
the entries of $s_\A(S)$, using tensoring over $\Lc$, and pull them
back to the left through the entries of $K_2$. This produces the desired effect for
right $t_\Bg(\tilde S_{op})$. We do the same with the  left matrix $t_\Bg(\tilde S_{op})$.
This time it is a bit more complicated, because we have to pass through $s_\A(S)$ obtained in the
previous step and through another matrix $K_2$. The result will be the same: the left matrix $t_\Bg(\tilde S_{op})$
is replaced by $s_\A(\tilde S)$.
Thus (\ref{eq_aux1}) transforms into
\be
\bar T_2 \bar T_1 s_\A(\tilde S) K_2 s_\A(\tilde S)  K_2 T_1 T_2
=
\bar T_2 \bar T_1  K_2 s_\A(\tilde S) K_2 s_\A(\tilde S) T_1 T_2,
\nn
\ee
where we used (\ref{dRE}).
Next we act in the reversed direction and replace $s_\A(\tilde S)$ by
$t_\Bg(\tilde S_{op})$. Then (\ref{eq_aux1}) becomes
\be
\bar T_2 \bar T_1  K_2 t_\Bg(\tilde S_{op}) K_2 t_\Bg(\tilde S_{op}) T_1 T_2
=\bar T_2 K_2 \bar T_1 t_\Bg(\tilde S_{op}) T_1 K_2  T_2 s_\Bg(\tilde S).
\nn
\ee
This is equal to $(\bar T_2 K_2 T_2)  s_\Bg(\tilde S)  (\bar  T_2 K_2  T_2 )s_\Bg(\tilde S)$,
so the relations (\ref{dRE}) are preserved.

We have shown that $\delta_{\tilde \K}$ is an algebra homomorphism.
To finish the proof, we must check conditions
2 and 3 of Definition \ref{comodule}. This is an easy and straightforward exercise.
\end{proof}
\section{Twisting the FRT and RE algebras}
\label{secTwFRTREa}
The remainder of the paper is devoted to the study of the dFRT and dRE algebras when the matrix $\tilde S$ is
obtained from a universal dynamical R-matrix.  That R-matrix itself is assumed to be a dynamical twist
of a quasitriangular structure on a Hopf algebra, $\U$. Our gual is to show that
the dFRT and dRE algebras are related to the ordinary FRT and RE algebras by certain dynamical twists.
As to the dFRT algebra, this is more or less straightforward and can be readily extracted
from \cite{DM2}. The case of dRE algebra is not obvious, so this will be in the focus of
our attention.
For better understanding of what is going on in the dynamical situation we first take a look
at how a Hopf algebra twist transforms the ordinary FRT and RE algebras.
\subsection{Twist  and twisted tensor square}
Let $\U$ be a quasitriangular Hopf algebra with the universal R-matrix
$\Ru$. Let $\U^{op}$ denote the coopposite Hopf algebra.
The FRT algebra is a module over $\U^{op}\tp \U$. Recall from \cite{DM3} that the RE algebra
a module over \tw{\U}{\Ru}{\U}, the twisted tensor square of $\U$, \cite{RS}.
Those two Hopf algebras  are related by the composition
of twists:
$$\U^{op}\tp \U \stackrel{\Ru_{13}}{\longrightarrow}\U\tp \U \stackrel{\Ru_{23}}{\longrightarrow} \mbox{\tw{\U}{\Ru}{\U}},$$
where $\Ru_{ij}\in (\U\tp \U)\tp (\U\tp \U)$.
This composite twist transforms the FRT algebra to the RE one.
So the latter is a twist of a module algebra over the intermediate Hopf algebra $\U\tp \U$.

Now suppose $\F$ is a twisting cocycle in $\U$, i.e. an invertible element from $\U\tp \U$
satisfying
\be
\F_{(12)3}\F_{12}=\F_{1(23)}\F_{23},
\quad
\ve_1(\F)=\ve_2(\F)=1_{\U^{\tp 2}}.
\label{tw_c}
\ee
Two twists $\F^{\{i\}}$, $i=1,2$, are called gauge equivalent if there is an invertible
element $v\in \U$ such that
\be
\Delta(v)\F^{\{1\}}=\F^{\{2\}}(v\tp v).
\label{gauge0}
\ee
If  $\tilde \U^{\{i\}}$, $i=1,2$,  be  respectively the twists of $\U$ by $\tilde \F^{\{i\}}$, then
the conjugation $x\mapsto v x \bar v$ implements a Hopf algebra isomorphism $\tilde \U^{\{1\}}\to \tilde \U^{\{2\}}$.

Given a twisting cocycle $\F$ of $\U$, we can consider the following
two twists $\Psi^{\{i\}}$, $i=1,2$, of the Hopf algebra $\U\tp \U$.
First we take the twist by $\Ru_{23}\in (\U\tp \U)\tp (\U\tp \U)$; then
the coproduct $\Delta\colon \U\to \mbox{\tw{\U}{\Ru}{\U}}$ is a Hopf algebra map and sends
the $\U$-cocycle $\F$ to the \tw{\U}{\Ru}{\U}-cocycle $\F_{(12)(34)}$.
We put $\Psi^{\{1\}}:=\F_{13}\F_{24}\tilde\Ru_{23}$.
The other twist is the composition of the twists $\F\tp \F$ and $\tilde \Ru_{23}=\F_{32}\Ru_{23}\F_{23}$,
thus we put $\Psi^{\{2\}}:=\Ru_{23}\F_{(12)(34)}$.

The question is how $\Psi^{\{1\}}$ and $\Psi^{\{2\}}$  are related to each other.
The answer is simple but very important for further considerations.
\begin{propn}
\label{coboundary}
The twists $\Psi^{\{1\}}$ and $\Psi^{\{2\}}$ are gauge equivalent:
\be
\F_{(13)(24)}\Psi^{\{1\}}=\Psi^{\{2\}}\F_{12}\F_{34}.
\label{eq_coboundary}
\ee
\end{propn}
\begin{proof}
Formula (\ref{eq_coboundary}) is a specialization of
the formula (\ref{gauge0}) for $\U\tp \U$ instead of $\U$,
with $v=\F$.
By definition of $\Psi^{\{1\}}$ and $\Psi^{\{2\}}$, equation (\ref{eq_coboundary})
is nothing else than
$$
\F_{(13)(24)}\F_{13}\F_{24}\tilde\Ru_{23}=\Ru_{23}\F_{(12)(34)}\F_{12}\F_{34}.
$$
But this latter equation is a corollary of the identity
$\F_{(12)(34)}\F_{12}\F_{34}=\F_{1(234)}\F_{2(34)}\F_{34}=\F_{(123)4}\F_{(12)3}\F_{12}$
following from the twisting cocycle equation (\ref{tw_c}).
\end{proof}
\subsection{How a twist of the Hopf algebra affects the RE algebra}
We are going to obtain the dRE algebra by a dynamical twist of the ordinary RE algebra.
To develop a strategy of solving this problem, we take a closer look at how the RE algebra behaves
under a non-dynamical or ordinary Hopf algebra twist.
The following diagram displays the relations between the FRT and RE algebras
associated with a quasitriangular Hopf algebra $(\U,\Ru)$ and their transformations under the
twist $\F$ of $\U$.
\be
\dgARROWLENGTH=0.8\dgARROWLENGTH
\begin{diagram}
\node{\Tg,\quad\U_{op}\tp \U}\arrow[2]{s,t}{\F^{-1}\tp \F}
\node{\U^{op}\tp \U}\arrow{w,t}{\gm\tp \id}\arrow[2]{s,l}{\F_{21}\tp \F}\arrow{e,t}{\Ru\tp \id}
\node{\U\tp \U}\arrow{e,l}{\Ru_{23}}\arrow[2]{s,t}{\F\tp \F}
\node{\mbox{\tw{\U}{\Ru}{\U}},}\arrow{se,t}{\Delta(\F)}
\node{\hspace{-1.5in}\Kg}
\\
\node[5]{\widetilde{\mbox{\tw{\U}{\Ru}{\U}}},\quad\tilde \Kg'}\arrow{sw,b}{\F^{-1}(\;)\F}
\\
\node{\tilde\Tg,\quad\tilde\U_{op}\tp \tilde\U}
\node{\tilde\U^{op}\tp \tilde\U}\arrow{w,t}{\tilde \gm\tp \id}\arrow{e,t}{\tilde \Ru\tp \id}
\node{\tilde\U\tp \tilde\U}\arrow{e,t}{\tilde\Ru_{23}}
\node{\mbox{\tw{\tilde\U}{\tilde\Ru}{\tilde\U}},}
\node{\hspace{-1.5in}\tilde\Kg}
\end{diagram}
\nn
\ee
The rightmost vertex of this diagram exists due to Proposition \ref{coboundary}.
Its outgoing arrow denotes the coboundary twist with $\F$ considered as an element of
\tw{\U}{\Ru}{\U}

In case of dynamical twist $\F$, we cannot draw exactly the same diagram.
We can construct dFRT algebra as a "bimodule"
over certain quantum groupoid extending $\U$ (see Subsection \ref{ssDFRT-T}), but
there are two severe obstructions to further proceeding to the dRE algebra.
First, the notion of bialgebroid antipode is not obvious.
Second, and this is crucial, there is no natural construction
of the twisted tensor square for quasitriangular bialgebroids.
Thus we cannot reach the right bottom corner moving counterclockwise from the left upper one.

Within the Hopf algebra setting, there is a way of getting $\tilde \Kg$ directly from $\Kg$
by applying two consecutive twists. Namely,
starting from the right upper corner and moving through the
right-most vertex.
In this passage, essential is the first twist, because the second one
is coboundary and results in an isomorphism of the module algebras.
The first or essential twist may well be constructed out of a dynamical cocycle $\F$
thus yielding an algebra, $\tilde \Kg'$, that may be called a dynamical RE algebra.
Our goal is to prove that this algebra is a specific case of the
dRE algebra introduced in Subsection \ref{ssDFRTa}.
To solve this problem, we should find a set of generators in
$\tilde \Kg'$ satisfying relations (\ref{Ktrans}) and
(\ref{dRE}) with $\tilde S$ being the image of the universal dynamical
R-matrix multiplied by the matrix permutation.

Let us consider in more detail how the RE transforms under the ordinary twist.
The RE algebra is commutative with respect to the universal R-matrix of \tw{\U}{\Ru}{\U},
and the RE relations are the corollary of this fact.
The universal R-matrix of \tw{\U}{\Ru}{\U} is expressed through
$\Ru$ and $\Ru^-:=\Ru_{21}^{-1}$ by $\Ru^-_{12}\Ru_{24}\Ru^-_{13}\Ru_{23}$.
The transition from
\tw{\U}{\Ru}{\U} to $\widetilde{\mbox{\tw{\U}{\Ru}{\U}}}$
destroys the RE relations in $\tilde\Kg'$.
 One reason for that is the
R-matrix of $\widetilde{\mbox{\tw{\U}{\Ru}{\U}}}$ loosing
its factorized form under this twist.
This problem can be fixed by the subsequent coboundary
twist from Proposition \ref{coboundary}. This transformation
restores the desired factorized form $\tilde \Ru^-_{12}\tilde \Ru_{24}\tilde \Ru^-_{13}\tilde \Ru_{23}$
of the universal R-matrix.
However this is not sufficient.
While \tw{\U}{\Ru}{\U} and \tw{\tilde\U}{\tilde\Ru}{\tilde\U} coincide
with $\U\tp \U$ as associative algebras, $\Kg$ and $\tilde\Kg$
are different (although isomorphic) as $\U\tp \U$-modules. Indeed, the action involves antipode, which
is changed under
the twist by a similarity transformation. Recall that the new antipode have the form
$\tilde \gm(h)=\gm(\bar\zt h\zt)$, where $\zt=\bar \F_2\bar\gm(\bar \F_1)$, \cite{Dr3}.
Therefore, if we perform the coboundary twist of \tw{\tilde\U}{\tilde\Ru}{\tilde\U}
induced by $\zt$, we bring the $\U\tp \U$-module structure of  $\tilde\Kg$
right to that of $\Kg$.

Summarizing, in the  Hopf algebra  case, the algebra  $\tilde\Kg$ can be obtained from $\Kg$ by
the twist
$(\Delta)(\F)$ of the Hopf algebra \tw{\U}{\Ru}{\U} and the subsequent coboundary twist induced by
$\F(\zt\tp \id)\in \mbox{\tw{\U}{\Ru}{\U}}$.
The first twist naturally carries over to the dynamical situation. It turns out
that the coboundary  twist also has a dynamical version. It will be constructed
in Sections \ref{secGTDT} and \ref{secDTREa}.

\section{Algebra of symmetric dynamical tensors}
\label{secASDT}
The present section is of technical character.
Here we study subspaces in  $\U^{\tp n}\tp \Lc$ possessing certain symmetries
with respect to the Hopf algebra $\Hg\subset \U$.
\subsection{Some definitions and notation}
\label{subsecSDN}
In the sequel of the paper we fix a quasitriangular
Hopf algebra $\U\supset \Hg$ with the universal R-matrix $\Ru$.
We will pursue a thorough study of a dynamical twist $\F$ over an $\Ha$-base algebra
$\Lc$ with
values in $\U$ and its certain derived twists.
When applied to the ordinary FRT and RE algebras
related to $\U$, those twists will be shown to produce the dFRT an dRE algebras
introduced in Section \ref{secdFRTREa}.
To this end, we need to develop a certain algebraic machinery.

Let us introduce notation $\check \U:=\U\tp \Lc$ and $\check \U^n:=\U^{\tp n}\tp \Lc$ for $n=2,3\ldots$.
It is convenient to put $\check \U^0:=\Lc$.
Elements of $\check \U^n$ will be called \select{dynamical tensors} of rank $n$.
We will use the following convention: $\U$-factors of $\check \U^n$ are labelled by positive integers;
the $\Lc$-factor is not indicated explicitly, assumed to be always on the right.
$\U$-factors stemming from the $\Ha$-coaction $\delta$ with the subsequent embedding to
$\U$ are marked by the right group of subscripts
separated by $|$, for example
$$
u_{1|2}:=u_1\tp u_\Lcc^{(1)}\tp u_\Lcc^{[\infty]},
\quad
u_{(12)|3}:=u^{(1)}_1\tp u^{(2)}_1\tp u_\Lcc^{(1)}\tp u_\Lcc^{[\infty]},
\quad
u_{31|42}:=u_2\tp u_\Lcc^{(2)}\tp u_1 \tp u_\Lcc^{(1)}\tp u_\Lcc^{[\infty]}
.
$$

Let $\S_n$ be the symmetric group of permutations of $n$-tuples. It
naturally acts on $\check \U^n$ by permutations of tensor factors. Denote by $\e$ the
unit element of $\S_n$ and by $\tau$ the flip $(1\ldots n)\mapsto (n\ldots 1)$.
For any pair $\al,\bt\in \S_n$, select in $\check \U^n$ the following subspace of $\Hg$-invariant elements:
$$
\Us{\al}{\bt}{n}:=\bigl\{u_{1\ldots n}\tp u_\Lcc \in\check\U^n\>|\>
(\Delta^n_\al)(h^{(1)})u_{1\ldots n}\tp h^{(2)}\tr u_\Lcc=
u_{1\ldots n}(\Delta^n_\bt)(h)\tp u_\Lcc, \forall h\in \Hg\bigr\}.
$$
Here $\Delta^n_\al:=\al\circ\Delta^n$, $\al\in \S_n$, and $\Delta^n$ is the $n$-fold coproduct.

We will say that a dynamical tensor of rank $n$ has \select{definite type of
symmetry} or just symmetric dynamical tensor if it belongs to $\Us{\al}{\bt}{n}$ for some $(\al,\bt) \in \S_n\times \S_n$.
There is a partial associative operation on the subspace of symmetric dynamical tensors, which
restricts from the multiplication in $\check \U^n$. We have,
 $\Us{\al}{\bt}{n} \Us{\bt}{\si}{n}\subset \Us{\al}{\si}{n}$ for
all $\al,\bt,\si\in \S_n$.
The $\S_n\times \S_n$-graded linear space $\check \U^n_{sym}=\oplus_{\al,\bt\in \S_n} \Us{\al}{\bt}{n}$ together
with this partial multiplication is called \select{algebra of symmetric dynamical tensors} of rank $n$.
Clearly $\Us{\al}{\al}{n}$ for any $\al\in \S_n$ is a unital associative algebra in the usual
sense. We consider the space  $\check \U^n$ of all dynamical tensors as a "two-sided module" over
$\check \U^n_{sym}$.

It is possible to introduce a structure of associative algebra on the infinite
direct sum $\oplus_{m=0}^\infty\check \U^m_{sym}$
which is a kind of tensor product. We will mostly use only its restriction to
$\oplus_{m=0}^\infty\>\Us{\e}{\e}{m}$.

\begin{propn}
\label{t_prod}
Let $i,k,n$ be positive integers such that  $i+k\leq n$. Then the embedding
\be
\label{embedding}
\iota_i\colon\check \U^k\to \check \U^{n},\quad \iota_i\colon u_{1\ldots k}\mapsto u_{i\ldots i+k-1|i+k\ldots n}
\ee
induces an embedding $\Us{\e}{\e}{k}\to \Us{\e}{\e}{n}$.
The map
\be
\label{t_prod_form}
\check \U^i\tp \check \U^k \to \check \U^{i+k},
\quad u\tp v \mapsto \iota_1(u)\iota_{i+1}(v)=\iota_{i+1}(v)\iota_1(u)
,\quad i,k=0,1\ldots
\ee
defines an associative
operation (tensor product) on $\oplus_{m=0}^{\infty}\>\Us{\e}{\e}{m}$.
\end{propn}
\begin{proof}
The embedding (\ref{embedding}) is obvious.
Due to (\ref{q-com}), elements from $\Us{\e}{\e}{k}$ commute with $\delta^k (\Lc)$.
Hence the image $\iota_1(\check \U^i)$ commutes with $\iota_{i+1}(\Us{\e}{\e}{k})$ in (\ref{t_prod_form}).
Associativity of the operation (\ref{t_prod_form}) follows from associativity and coassociativity of $\delta$.
\end{proof}
\begin{remark}
Note that symmetric dynamical tensors of type $(\e,\e)$ yield morphisms in dynamical categories
of \cite{DM1}. The tensor product from Proposition \ref{t_prod} is induced by the tensor
product of morphisms.
\end{remark}
\subsection{Dynamical twist and R-matrix}
Let $\U$ be a Hopf algebra containing $\Hg$ as a Hopf subalgebra.
A dynamical $(\U,\Hg,\Lc)$-twist is an element $\F\in \Us{\e}{\e}{2}$ satisfying the
conditions
\be
\label{sh_cocycle}
&\F_{(12)3}\F_{12|3}=\F_{1(23)}\F_{23},
\\
&\ve_1(\F)=1_{\check \U^2} =\ve_2(\F).
\label{norm}
\ee
Numerous examples of dynamical twists over various bases relative both to classical and quantum groups are built in
\cite{EV1,ESS,EE,EEM}.

The notion of coboundary twist can be modified for the dynamical situation. Specifically,
two twists $\F^{\{i\}}\in \Us{\e}{\e}{2}$, $i=1,2$,
are called gauge equivalent if there exists $\zeta\in \Us{\e}{\e}{}$ such that
$$
\Delta(\zeta)\F^{\{1\}}=\F^{\{2\}}\zt_{1|2}\zt_2.
$$
Note that each factor in this equality has the symmetry type $(\e,\e)$;
so this is an equality in the associative algebra $\Us{\e}{\e}{2}$. Also note that
 two $\zt$-factors on the right commute with each other, see (\ref{t_prod_form}).

Observe that  the R-matrix $\Ru$ of the Hopf algebra  $\U$ belongs to $\Us{\tau}{\e}{2}$.
The element $\bar\F_{21}\Ru\F\in \Us{\tau}{\e}{2}$ is called dynamical R-matrix and satisfies the
(dynamical Yang-Baxter) equation
\be
\tilde \Ru_{12} \;\tilde\Ru_{13|2}\; \tilde \Ru_{23}&=&
\tilde \Ru_{23|1}\; \tilde \Ru_{13}\;  \tilde \Ru_{12|3}.
\label{uR}
\ee
Note that each factor in (\ref{uR}) has a definite symmetry type in $\check \U^3_{sym}$,
thus (\ref{uR}) is an equation in $\check \U^3_{sym}$.

The gauge transformation of two dynamical twist corresponding to the element $\zt\in \Us{\e}{\e}{}$ amounts to the
transformation of a dynamical R-matrix
$$
\tilde \Ru\mapsto \bar \zt_{2|1}\bar\zt_1 \tilde \Ru\zt_{1|2}\zt_2.
$$
The dynamical R-matrix $\tilde \Ru$ defines  braiding in certain monoidal categories, \cite{DM1,DM2}.

\subsection{One technical lemma}
Denote $\omega=\Theta^{-1}_1\gm^{-1}(\Theta^{-1}_2)$, a Drinfeld element of the double
Hopf algebra $\D\Hg$, \cite{Dr2}.
It satisfies the equality $\Delta(\omega^{-1})(\omega\tp \omega)=\Theta^{-1}\Theta^{-1}_{21}$.
From this we conclude that the map $\la\mapsto \omega\tr \la$ is an automorphism of the base algebra $\Lc$
(recall that  in terms of $\Theta$ the  condition (\ref{q-com}) reads $(\Theta_2\tr\la)(\Theta_1\tr\mu)=\mu\la$ for all $\la,\mu\in \Lc$).
Conjugation with $\omega$ implements the squared antipode in the Hopf algebra $\D\Hg$.

Let $\A$ be an associative algebra and let $(\la,\rho)$ be a pair of homomorphisms
from $\Hg$ to $\A$.
Denote by $(\A\tp \Lc)_{\ell,\rho}$ the subset of elements $a=a_1\tp a_\Lcc\in \A\tp \Lc$ satisfying
\be
\label{l-r}
\ell(h^{(1)})a_1\tp h^{(2)}\tr a_\Lcc=a_1 \rho(h)\tp a_\Lcc, \quad \forall h\in \Ha.
\ee
Equivalently, condition (\ref{l-r})  can be represented as
\be
\label{l-r0}
\ell(h)a_1\tp a_\Lcc=a_1 \rho(h^{(1)})\tp \bar\gm(h^{(2)})\tr a_\Lcc, \quad h\in \Hg.
\ee
It follows that
\be
\label{l-r1}
a_1 \rho(\Theta_2)\tp \bar\gm(\Theta_1)\tr a_\Lcc= \ell(\bar\Theta_2) a_1\tp (\bar\omega\bar\Theta_1)\tr a_\Lcc
\ee
for all
$a=a_1\tp a_\Lcc\in (\A\tp \Lc)_{\ell,\rho}$. Indeed, we have for the right-hand side of
(\ref{l-r1})
\be
\ell(\bar\Theta_2) a_1\tp (\bar\omega\bar\Theta_1)\tr a_\Lcc
&=&
a_1\rho(\bar\Theta^{(1)}_2) \tp \bigl(\bar\omega\bar\Theta_1  \bar\gm(\bar\Theta^{(2)}_2)\bigr)\tr a_\Lcc
\nn\\
&=&
a_1\rho(\bar\Theta_{2}) \tp \bigl(\bar\omega\bar\Theta_{1}\bar\Theta_{1'}\bar\gm(\bar\Theta_{2'})\bigr)\tr a_\Lcc
=
a_1\rho(\bar\Theta_{2}) \tp (\bar\omega\bar\Theta_{1}\omega)\tr a_\Lcc.
\nn
\ee
In the first equality of this chain we used (\ref{l-r0}).
It remains to resort to the properties of $\omega$, $\gm$, and $\Theta$ and get
the left-hand side of(\ref{l-r1}).
\begin{lemma}
\label{lm_auxi0001}
For all $a=a_1\tp  a_\Lcc\in (\A\tp \Lc)_{\ell,\rho}$ and
for any  $\mu\in \Lc$:
\be
\ell(\bar \Theta_2) a_1 \tp \bar \Theta_1\tr( a_\Lcc\mu)=
\ell(\bar \Theta_{2'}) a_1 \rho(\bar \Theta_2) \tp (\bar \Theta_1\tr\mu)(\bar \Theta_{1'}  \tr a_\Lcc)
.
\label{eq_aux0001}
\ee
\end{lemma}
\begin{proof}
Condition (\ref{eq_aux0001}) is equivalent to
\be
\ell(\bar \Theta_2) a_1 \tp (\bar \omega\bar \Theta_1)\tr( a_\Lcc\mu)=
\ell(\bar \Theta_{2'}) a_1 \rho(\bar \Theta_2) \tp \bigl((\bar \omega\bar \Theta_1)\tr\mu\bigr)
\bigl((\bar \omega\bar \Theta_{1'})  \tr a_\Lcc\bigr )
\label{eq_aux0}
\ee
because $\omega$ implements a homomorphism of $\Lc$. So let us prove  (\ref{eq_aux0}) instead of (\ref{eq_aux0001}).

The proof  employs the standard machinery of quasitriangular Hopf algebras as specialized to
the double $\D\Hg$.
In view of  (\ref{l-r0}) and (\ref{l-r1}), the left-hand side of (\ref{eq_aux0}) can be transformed as
\be
&&\ell(\bar \Theta_{2} \bar \Theta_{2'}) a_1 \tp \bigl((\bar \omega\bar \Theta_{1'})\tr  a_\Lcc\bigr)
\bigl((\bar \omega\bar\Theta_{1})\tr\mu\bigr)
=
\ell(\bar \Theta_{2}) a_1  \rho(\Theta_{2'})  \tp
\bigl((\bar\gm(\Theta_{1'})\tr  a_\Lcc\bigr)\bigl((\bar \omega\bar\Theta_{1})\tr\mu\bigr)
\nn\\
&&\hspace{2cm}=
a_1 \rho(\bar \Theta_{2}^{(1)}  \Theta_{2'})  \tp
\bigl((\bar\gm(\Theta_{1'})\bar\gm(\bar \Theta_{2}^{(2)})\tr  a_\Lcc\bigr)\bigl((\bar \omega\bar\Theta_{1})\tr\mu\bigr)
\nn\\
&&\hspace{2cm}=
a_1 \rho(\Theta_{2'} \bar \Theta_{2}) \tp \bigl(\bar \gm(\bar \Theta_{2''})\bar\gm(\Theta_{1'})\tr  a_\Lcc\bigr)
\bigl((\bar \omega\bar\Theta_{1''}\bar\Theta_{1})\tr\mu\bigr).
\nn
\label{eq_aux01}
\ee
The last equality is the result of permutation of $\Theta$ with the coproduct
and the subsequent factorization of $(\id\tp\Delta)(\Theta)$.
Now we pull $\bar \omega$ through $\bar\Theta_{1''}$ to the right and use the $\gm$-symmetry of $\Theta$.
After that we employ $\Theta$-commutativity of the base algebra $\Lc$ and bring (\ref{eq_aux01}) to
$a_1 \rho(\Theta_{2'} \bar \Theta_{2}) \tp \bigl((\bar \omega\bar\Theta_{1})\tr\mu \bigr)\bigl(\gm^{-1}(\Theta_{1'})\tr  a_\Lcc\bigr)$.
Now, to get the right-hand side of (\ref{eq_aux0}), we again apply (\ref{l-r1}) in the reversed direction.
\end{proof}

\subsection{An anti-automorphism of $\check{\mathcal{U}}^n_{sym}$}
In this subsection we introduce an important family of endomorphisms of dynamical tensors.
This family yields a kind of  anti-automorphism of the "algebra" of symmetric dynamical tensors.

Introduce a linear operator $\Gm\in \Aut(\check \U)$ setting
$$
\Gm (x\tp \la)
:=\gm (\bar \Theta_2 x)\tp \bar\Theta_1\tr \la
$$
for all $x\tp \la\in \check \U$. Further, define a linear automorphism  $\Gm_i$ of $\check \U^n$ as $\Gm$ acting
on the i-th $\U$ and the $\Lc$-factor of $\check \U^n$.
The operators $\Gm_i$ and $\Gm_j$ do not commute
unless $i =j$.
Define a linear automorphism $\Gm_{1\ldots n}\in \Aut(\check \U^n)$ for $n\in \N$ setting $\Gm_{1\ldots n}:=\Gm_1\ldots \Gm_n$.
Put also $\Gm_{\si}:=\si\circ\Gm_{1\ldots n}$ for $\si\in \S_n$; then
the operator $\Gm_{\si}$ can be represented as
\be
\Gm_{\si}(x_{1\ldots n}\tp x_\Lcc)=\gm^{\tp n}\bigl(\Delta^n_\si(\bar\Theta_2)x_{1\ldots n}\bigr)\tp \bar\Theta_1\tr x_\Lcc
\quad
\mbox{for }
x=x_{1\ldots n}\tp x_\Lcc \in \U^{\tp n}\tp \Lc =\check \U^n.
\label{Gm_rep}
\ee
Similarly we define the operators $\bar\Gm_{1\ldots n}:=\bar\Gm_1\ldots \bar\Gm_n$ and
$\bar \Gm_{\si}:=\si\circ\bar\Gm_{1\ldots n}$ for $\si\in \S_n$.
We have $(\Gm_\si)^{-1}=\bar \Gm_{\si'}$, where
the map $\al\mapsto \al'$ is the inner automorphism of $\S_n$ with the element $\tau$,
namely
$\al'=\tau \al \bar\tau$.

The collection $\{\Gm_\al\}_{\al\in \S_n}$ is an anti-homomorphism of the algebra
of symmetric tensors. Moreover, it implements an isomorphism between
the left and right $\check \U^n_{sym}$-"module structure" on $\check \U^n$.
To be precise, we have the following proposition.
\begin{propn}
\label{Gm}
(i) Fix a pair of  permutations $\al,\bt\in \S_n$. Then
\be
\Gm_\al\bigl(\Us{\al}{\bt}{n}\bigr)\subset \Us{\bt'}{\al'}{n},
\quad
\bar\Gm_\al\bigl(\Us{\al}{\bt}{n}\bigr)\subset \Us{\bt'}{\al'}{n}
\quad \mbox{and}
\label{Gm001}
\\
\Gm_\al(x y)=\Gm_\bt(y)\Gm_\al(x)
,
\quad
\bar \Gm_\al(y x)=\bar \Gm_\al(x) \bar\Gm_\bt(y)
\label{Gm002}
\ee
for
all $x\in \Us{\al}{\bt}{n}$ and $y\in \check \U^n$.
\\
(ii) For all $x\in \check \U$ one has
\be
(\Gm_\si\circ\Delta^n_\si)(x)=(\Delta^n_{\si'}\circ \Gm)(x),
&&
(\bar\Gm_\si\circ\Delta^n_\si)(x)=(\Delta^n_{\si'}\circ \bar\Gm)(x),
\label{Gm003}
\ee
\end{propn}
\begin{proof}
The right formulas in (\ref{Gm001}-\ref{Gm003}) follow from the left ones, as
$(\Gm_\si)^{-1}=\bar \Gm_{\si'}$. Let us sketch the proof for the left ones.
The left inclusion (\ref{Gm001}) is a corollary of
(\ref{Gm_rep}) and $(\Gm_\si)^{-1}=\bar \Gm_{\si'}$.
The left equality (\ref{Gm002}) follows from Lemma  \ref{lm_auxi0001}. These prove (i).
Statement (ii) follows from (\ref{Gm_rep}).
The details are left to the reader.
\end{proof}
\begin{example}
\label{autdYBE}
As an application of Proposition \ref{Gm}, apply $\Gm_{321}$ to the dYBE  (\ref{uR}),
which is an equation in $\check U^3_{sym}$.
It is easy to see that $\Gm_{321}$ implements an automorphism of
(\ref{uR}), and the correspondence $\tilde \Ru\mapsto\Gm_{21}(\tilde \Ru)$
maps a solution to another solution.
\end{example}
\section{A gauge transformation of dynamical twist}
\label{secGTDT}
In the theory of Hopf algebras, a twist transformation
changes the comultiplication and the antipode by a similarity
transformation. The elements $\vt$ and $\bar \vt=\vt^{-1}$ defining the new antipode
$\tilde \gm(h)= \bar\vt^{-1} \gm(h) \vt$ are expressed
through the  twisting cocycle $\F$  as $\vt=\gm(\F_1)\F_2$ and $\bar\vt=\bar\F_1\gm(\bar\F_2)$.
They implement the gauge-equivalence between the twists $\F$ and $(\gm\tp\gm)(\bar\F_{21})$,
see \cite{Dr3}:
\be
\label{gauge}
\Delta(\vartheta)\F = (\gm\tp\gm)(\bar\F_{21})(\vartheta\tp \vartheta).
\ee
It turns out that the dynamical version of the element $\vt$,
as well as the equation (\ref{gauge}), does exist. They play an important role in our
consideration.

\subsection{Elements $\vt$ and $\zt$}
Let $\F=\F_1\tp\F_2\tp\F_{\Lcc}\in \Us{\e}{\e}{2}$ be a dynamical twisting cocycle.
Introduce the elements $\vt,\zt\in \U\tp \Lc$ by
\be
\vartheta:=\gm(\F_1)\F_2\tp \F_\Lcc=\vartheta_1\tp \vartheta_\Lcc
,\quad
\zt:=\bar \F_{2}\bar\gm(\bar \F_{1})\tp  \bar \F_{\Lcc}
=\zeta_1\tp \zeta_\Lcc.
\label{tz}
\ee
Since  $\F\in \Us{\e}{\e}{2}$, one can show that $\vt$ and $\zt$ belong to $\Us{\e}{\e}{}$.
\begin{lemma}
Let $\F$ be a dynamical twisting cocycle and $\vt$, $\zt$ defined  by
(\ref{tz}). Then
\be
\F^{(1)}_1\tp \gm(\F^{(2)}_1)\F_2\tp \F_{\Lcc}
&=&
\Gm_2(\bar\F)\vt_2,
\label{cc1}\\
\gm(\F_1)\F^{(1)}_2\tp \F^{(2)}_2\tp  \F_{\Lcc}
&=&
\vt_{1|2}
\bar \F,
\label{cc2}\\
\bar \F^{(1)}_2\bar \gm(\bar \F_1)\tp \bar \F^{(2)}_2\tp  \bar \F_{\Lcc}
&=&
\F
\zt_{1|2},
\label{cc3}\\
\bar\F^{(1)}_1\tp \bar\F_2\bar \gm(\bar\F^{(2)}_1)\tp \bar\F_{\Lcc}
&=&
\zt_2
\bar\Gm_2(\F).
\label{cc4}
\ee
\end{lemma}
\begin{proof}
Direct consequence of the conditions (\ref{sh_cocycle}) and (\ref{norm}).
\end{proof}
\begin{propn}
\label{prop_vt_inv}
The elements $\vt$ and $\zt$ are invertible and
\be
\bar\vt
\label{vt_inv}
=
\Gm (\zt),
\quad
\bar\zt=\bar\Gm (\vt).
\ee
\end{propn}
\begin{proof}
Note that  one of these equations implies the other since
$\Gm$ is an anti-automorphism of $\Us{\e}{\e}{}$.

We will deduce the left equality (\ref{vt_inv}) from (\ref{cc1}) and (\ref{cc4}).
First observe that taking product of the $\U$-factors in $\Gm_2(\bar\F)$ gives
$\Gm(\zt)$. Then (\ref{cc1}) implies
$
1_{\check \U}=\Gm (\zt)\vt,
$
i.e.
$\Gm (\zt)$  is a left inverse to $\vt$.
From (\ref{cc4}) we obtain
$
1_{\check \U}=\zt \bar\Gm (\vt).
$
Apply to this equation the operator $\Gm$, which is an anti-algebra automorphism of $\Us{\e}{\e}{}$.
This shows that $\Gm (\zt)$ is also a right inverse to $\vt$.
\end{proof}
\subsection{The gauge transformation}
A dynamical twist induces an isomorphic transformation of the dynamical extension of the monoidal
category of $\U$-modules, \cite{DM1}. Therefore it can be composed with a
coboundary twist (a coboundary twist may be associated with an arbitrary
collection of automorphisms). We call such
coboundary twist a gauge transformation. A gauge transformation
of a dynamical twist $\F$ is defined by any invertible element $\xi\in \Us{\e}{\e}{}$
and reads
$$\F\mapsto \Delta(\bar\xi)\F \xi_{1|2}\xi_{2}.$$
\begin{propn}
The symmetric dynamical tensor $(\Gm_{21}) (\bar \F_{21})\in \Us{\e}{\e}{2}$ is gauge-equivalent to $\F$ with
\be
\Delta(\vartheta)\F
&=&
(\Gm_{21}) (\bar \F_{21})\>(\vartheta_{1|2}\vartheta_2),
\label{eq_v0}\\
\bar \F\Delta(\zt)
&=&
(\zt_{1|2}\zt_2)\>(\bar \Gm_{21})(\F_{21}).
\label{eq_v}
\ee
\end{propn}
\noindent
\begin{proof}
Let us compute the expression $\Delta(\vartheta)\F$.  Using
(\ref{sh_cocycle}), we find
$$
\Delta(\vartheta)\F=\Bigr(\Delta\bigl(\gm(\F_{1})\bigr)\Delta(\F_{2})\tp \F_\Lcc\Bigr)\F
=
\Delta\bigl(\gm(\F^{(1)}_{1}\F_{1'})\bigr)(\F^{(2)}_{1}\F_{2'}\tp \F_{2}\F^{(1)}_{\Lcc'} )\tp \F_\Lcc\F^{[\infty]}_{\Lcc'}
.
$$
Permuting the coproduct and antipode and using (\ref{norm}), we find this equal to
\be
\gm(\F^{(2)}_{1'})
\F_{2'}\tp
\gm(\F^{(1)}_{1'})
\gm(\F_{1})
\F_{2}\F^{(1)}_{\Lcc'} \tp \F_\Lcc\F^{[\infty]}_{\Lcc'}
=
\bigl(\gm(\F^{(2)}_{1'})\F_{2'}\tp
\gm(\F^{(1)}_{1'})\F^{(1)}_{\Lcc'} \tp \F^{[\infty]}_{\Lcc'}\bigr)\vt_2
\label{eq0001}.
\ee
Here we used the fact that $\vt$ belongs to $\Us{\e}{\e}{}$ and therefore commutes with $\delta(\F_{\Lcc'})$.
Now let us compute left factor in the right-hand side of (\ref{eq0001}).
From (\ref{cc1}) we deduce
\be
\gm(\F^{(2)}_1)\F_2\tp \gm(\F^{(1)}_1)\F^{(1)}_{\Lcc}\tp \F^{[\infty]}_{\Lcc}
&=&
\Gm_2(\Gm_1(\bar \F_{21})\vt_1)
=
\Gm_{21}(\bar \F_{21})
\vt_{1|2}.
\nn
\ee
The last equality can be checked directly.
Substituting this into (\ref{eq0001}), we prove the formula (\ref{eq_v0}).

Now we derive (\ref{eq_v}) from (\ref{eq_v0}).
First of all observe, using the right equality (\ref{vt_inv}) and Proposition \ref{Gm}, that
\be
(\bar\Gm_{12}\circ\Delta)(\vt)&=&(\Delta_{\tau}\circ \bar\Gm)(\vt)=
\Delta_{\tau}(\bar\zt),
\nn\\
(\bar\Gm_{12})(\vt_{1|2}\vt_2)&=&
(\bar\Gm_1\bar\Gm_2)(\vt_{2})(\bar\Gm_1\bar\Gm_2)(\vt_{1|2})=\bar \zt_{2|1}\bar \zt_1.
\nn
\ee
Note that $\Delta_\tau=\Delta_{op}$.
Here we employed the fact that $\vt_{1|2}$ and $\vt_2$ belong to $\Us{\e}{\e}{2}$.
Applying $\tau\circ\bar\Gm_{12}$ to (\ref{eq_v0}) and using the above equalities, we obtain
(\ref{eq_v}).
\end{proof}
\subsection{Symmetries of dynamical R-matrix}
In the theory of quasitriangular Hopf algebras,
the universal R-matrices transform under the action
of antipode as
$$
(\gm\tp\gm)(\Ru)=\Ru,
\quad
(\gm\tp\id)(\Ru)=\bar\Ru,
\quad
(\id\tp\gm)(\bar\Ru)=\Ru.
$$
 In the present subsection we establish dynamical analogs of these relations.
\begin{propn}
\label{Rmatsym}
Suppose that $\U$ is a quasitriangular Hopf algebra with the universal R-matrix $\Ru$.
Let $\F$ be a dynamical $(\U,\Hg,\Lc)$-twist, $\tilde \Ru=\F_{21}^{-1} \Ru\F$ the corresponding dynamical R-matrix, and
$\zeta$ defined as in (\ref{tz}). Then
\be
(\Gm_2\Gm_1)(\bar\zeta_{2|1}\bar\zeta_1\tilde \Ru\zeta_{1|2}\zeta_2)
=&
\tilde \Ru
&=
(\bar\Gm_1\bar\Gm_2)(\vt_{2|1}\vt_1\tilde \Ru\bar\vt_{1|2}\bar \vt_2)
,
\label{r_ant1}
\\
\Gm_1(\bar \zeta_1\tilde \Ru \zt_{1|2})
=&
\tilde \Ru^{-1}
&=
\bar\Gm_2(\vt_{2|1}\tilde \Ru \bar \vt_2),
\label{r_ant2}
\\
\Gm_2(\bar\zeta_2\tilde \Ru^{-1} \zt_{2|1})
=&
\tilde \Ru
&=
\bar \Gm_1(\vt_{1|2}\tilde \Ru^{-1} \bar\vt_1)
.
\label{r_ant3}
\ee
\end{propn}
\begin{proof}
Let us establish the left equality (\ref{r_ant2}) first. Using (\ref{cc3}) we find
\be
\bar\F_{21}\Ru\F\zeta_{1|2}
& =&
\bar\F_{2'}\Ru_1\bigl(\bar \F^{(1)}_2\bar \gm(\bar \F_1)\bigr)\tp \bar\F_{1'}\Ru_2(\bar \F^{(2)}_2)\tp\bar\F_{\Lcc'}(\bar \F_{\Lcc}).
\label{eq_aux2}
\ee
Pulling the term $\Delta(\F_2)$ through the R-matrix to the left and then
 applying the dynamical cocycle condition (\ref{sh_cocycle})
we transform the right-hand side of (\ref{eq_aux2}) to
$$
\bar \F^{(1)}_{\Lcc'}\bar\F_{2} \Ru_1 \bar \gm(\bar\F_{1'}\bar \F^{(1)}_1)\tp
\bar\F_{2'}\bar \F^{(2)}_1 \Ru_2\tp\bar\F^{[\infty]}_{\Lcc'}\bar \F_{\Lcc}
=
\bar \F^{(1)}_{\Lcc'}\bar\F_{2} \bar \gm(\bar \F^{(2)}_1) \Ru_1\bar \gm(\bar\F_{1'})\tp
\bar\F_{2'}\Ru_2\bar \F^{(1)}_1 \tp\bar\F^{[\infty]}_{\Lcc'}\bar \F_{\Lcc}
,$$
where we used the equality
$
\Ru_1 \bar \gm(\F^{(1)}_1)\tp \F^{(2)}_1 \Ru_2=\bar \gm(\F^{(2)}_1)\Ru_1 \tp \Ru_2\F^{(2)}_1
$.
Now we appeal to (\ref{cc4}) and find (\ref{eq_aux2}) equal to
$\bar \F^{(1)}_{\Lcc'}\zeta_1\F^{(1)}_{\Lcc''}\bar\gm(\F_{2''}) \Ru_1\bar \gm(\bar\F_{1'})\tp
\bar\F_{2'}\Ru_2\F_{1''} \tp\bar\F^{[\infty]}_{\Lcc'}\zeta_\Lcc\F^{[\infty]}_{\Lcc''}$.
Since $\zt\in \Us{\e}{\e}{}$, it commutes with $\delta(\Lc)$ and can be pulled to the left. Thus we come to
\be
\bar\F_{21}\Ru\F\zeta_{1|2}&=&\zeta_1\bar \F^{(1)}_{\Lcc'}\F^{(1)}_{\Lcc''}\bar\gm(\F_{2''}) \Ru_1\bar \gm(\bar\F_{1'})\tp
\bar\F_{2'}\Ru_2\F_{1''} \tp\zeta_\Lcc\bar\F^{[\infty]}_{\Lcc'}\F^{[\infty]}_{\Lcc''}
\nn\\
&=&
\zt_1\bigl(
\Theta_2\bar\gm(\tilde\Ru_1)\tp \tilde \Ru_2\tp \Theta_1\tr\tilde \Ru_\Lcc
\bigr)
.
\ee
This immediately implies the  left equality (\ref{r_ant2}).

Further, using (\ref{eq_v}), we find
$$
\bar\zeta_{2|1}\bar\zeta_1\tilde \Ru\zeta_{1|2}\zeta_2=
\bar \Gm_{12}(\F)\Ru \bar \Gm_{21}(\bar\F_{21}).
$$
Observe that each factor on the right-hand side has definite symmetry type,
and this product has the structure $\Us{\tau}{\tau}{2} \Us{\tau}{\e}{2} \Us{\e}{\e}{2}$.
Apply to this equation the map $\Gm_{21}=\Gm_2\Gm_1=(\bar\Gm_{12})^{-1}$
and use Proposition \ref{Gm} (i), formula (\ref{Gm002}). This proves the left equality (\ref{r_ant1}),
from which we deduce the right one. Indeed, the left equality (\ref{r_ant1}) implies
$$
\tilde \Ru= (\Gm_1\Gm_2)(\zeta_{1|2}\zeta_2)(\Gm_2\Gm_1)(\tilde \Ru) (\Gm_2\Gm_1)(\bar\zeta_{2|1}\bar\zeta_1).
$$
Here we again used Proposition \ref{Gm} (i), formula (\ref{Gm002}). It remains to
apply formulas (\ref{vt_inv}), to get the right equality (\ref{r_ant1}).

The right equality  (\ref{r_ant2}) follows from  (\ref{r_ant1}) and the left equality (\ref{r_ant2}).
Indeed, observe that the element $\bar\zeta_1\tilde \Ru\zeta_{1|2}$ belongs to $\Us{\tau}{\e}{2}$.
Employing (\ref{Gm002}) we find
$$
\tilde \Ru = (\Gm_2\Gm_1)(\bar\zeta_{2|1}\bar\zeta_1\tilde \Ru\zeta_{1|2}\zeta_2)
= (\Gm_1\Gm_2)(\zeta_2)(\Gm_2\Gm_1)(\bar\zeta_1\tilde \Ru\zeta_{1|2})(\Gm_2\Gm_1)(\bar\zeta_{2|1})
= \bar\vt_{2|1} \Gm_2 (\tilde\Ru^{-1})\vt_2.
$$
The last equality is obtained using (\ref{vt_inv}) and the left equality (\ref{r_ant2}).

The left and right equalities  (\ref{r_ant3}) are corollaries of the right and left equalities (\ref{r_ant2})
respectively. This completes the proof.
\end{proof}
\begin{remark}
\label{rigid}
The gauge equivalence (\ref{eq_v0})--(\ref{eq_v0}) as well as the formulas (\ref{r_ant1})--(\ref{r_ant3})
take place because of rigidity of the dynamical category of finite dimensional $\Ha$-modules.
Indeed, the dynamical category is rigid, being an extension of the rigid category of finite dimensional
 $\Ha$-modules by enlarging homsets, see \cite{DM1}.
\end{remark}

\section{Dynamical twist of the RE algebra}
\label{secDTREa}

\subsection{Dynamical twist and quantum groupoids}
\label{subsecdRQG}
All the bialgebroids in this subsection
are assumed to be left.
Here we recall the connection between the dynamical twist
and twist of bialgebroids following \cite{Xu1,DM2}.

Given a Hopf algebra inclusion $\Hg\subset \U$ and a dynamical $(\U,\Hg,\Lc)$-twist $\F$,
it is possible to construct a twisted left $\Lc$-bialgebroid $\widetilde{\U\tp \D\Hg_\Lc}$
starting from the tensor product $\Lc$-bialgebroid $\U\tp \D\Hg_\Lc$, see \cite{DM2}.
This latter bialgebroid is obtained  by a projection from $\U\tp \Lc\rtimes\D\Hg$, cf. Example \ref{ex_bialg3}.
We will use this projection representation for elements of  $\U\tp \D\Hg_\Lc$.
Under this convention, the bialgebroid twist ${\bf\Psi}\in (\U\tp \D\Hg_\Lc)\tp_\Lc(\U\tp \D\Hg_\Lc)$ is
expressed through $\F$ and $\Theta$ by
\be
\label{biatwist}
{\bf\Psi}&:=&(\F_1\tp \F_{3} \otimes \Theta_1)\tp_{\Lc}(\F_2\Theta_2\tp 1_\Lc\otimes 1_{\D\Hg}).
\ee
We assume that $\U$ is a quasitriangular Hopf algebra with the universal R-matrix $\Ru$.
Then the bialgebroids $\U\tp \D\Hg_\Lc$ and $\widetilde{\U\tp \D\Hg_\Lc}$ are also quasitriangular, \cite{DM2}.
The universal R-matrix of $\widetilde{\U\tp \D\Hg_\Lc}$ reads
\be
\label{R-matrix-twisted}
{\bf R}=
(\bar\Theta_{2'}\bar\Theta_{2'''}\tilde \Ru_1\tp 1_\Lc\otimes \Theta_1\Theta_{1''})
\tp_{\Lc_{op}}
(\tilde \Ru_2\Theta_{2''}\tp \bar\Theta_{1'''}\tr\tilde \Ru_3\otimes\bar\Theta_{1'}\Theta_2),
\ee
where $\tilde \Ru= \bar\F_{21} \Ru \F=\tilde \Ru_1\tp \tilde \Ru_2\tp \tilde \Ru_3\in \U\tp \U\tp \Lc$
is the dynamical R-matrix, the twist of $\Ru$.

Let $\A$ be a left $\U$-module algebra and $\btr$ denote the  $\U$-action on $\A$.
Recall that $\Lc$ is a module algebra over the bialgebroid $\D\Hg_\Lc$ with respect to
the anchor action.
Then the tensor product $\A\tp \Lc$ is a natural module (and module algebra)
over the tensor product bialgebroid $\U\tp \D\Hg_\Lc$.
Recall from Example \ref{ex_bialg3} that $\D\Ha$ is a quotient of $\Lc\rtimes \D\Ha$.
The action of $\U\tp \D\Hg_\Lc$ is induced by the action of $\U\tp \Lc\rtimes \D\Hg$,
which is given by
\be
(u\tp \la\tp h)\tp (a\tp \mu)\mapsto u\btr a\tp \la(h\tr \mu),
\ee
for $(u\tp \la\tp h)\in \U\tp \Lc\rtimes \D\Hg$ and $a\tp \mu\in\A\tp \Lc$.

Suppose that $\A$ is a commutative algebra in the braided category of $\U$-modules,
i.e.
$$
(\Ru_2\btr b)(\Ru_1\btr a)= ab,
\quad \mbox{for all} \quad a,b\in \A.
$$
Then
$\A\tp \Lc$ is a commutative algebra in the  braided category of $\U\tp \D\Hg_\Lc$-modules.
The bialgebroid twist defines an automorphism of the category of $\U\tp \D\Hg_\Lc$-modules,
so it transforms a commutative algebra into commutative, with respect to the new braiding.

The bialgebroid twist that comes from a dynamical one re-defines the multiplication
in the algebra $\A\tp \Lc$ according to  the rule
\be
\begin{array}{llllll}
a*b&=&(\F_1\btr a)(\F_2\btr b)\tp \F_\Lcc,
\quad
&\la*\mu&=&(1\tp \la\mu),
\\[5pt]
\la*a&=&\la^{(1)}\btr a\tp \la^{[\infty]},
\quad
&a*\la&=&a\tp \la.
\end{array}
\label{*prod}
\ee
for $a,b\in \A \subset \A\tp \Lc$ and $\la,\mu \in \Lc  \subset \A\tp \Lc$.
The new algebra $\tilde \A$ is that over the twisted quantum groupoid
$\widetilde{\U\tp \D\Hg_\Lc}$ and, in particular, it is an $\Lc$-bimodule algebra.
It is commutative with respect to the permutation induced by $\mathbf{R}$.

It is seen that the permutation rule between $\Lc\subset \tilde \A$ and the entire $\tilde \A$
is governed by $\Theta$, the R-matrix of the double $\D\Hg$.
Permutation among elements from $\A\subset \tilde \A$ is expressed through
the dynamical R-matrix by the following obvious formula.
\be
\label{dyn_alg}
(\tilde\Ru_2 \btr b)*(\tilde\Ru_1 \btr a)*\tilde\Ru_\Lcc*\la
=a*b*\la
\ee
for any $a,b \in \A$ and $\la\in \Lc$.
\subsection{dFRT algebra via dynamical twist}
\label{ssDFRT-T}
Let $(V,\rho)$ be a $\U$-module and therefore an $\Hg$-module by restriction.
Denote by $R$ and $\tilde R$ the images of $\Ru$ and $\tilde \Ru$ in this representation.
Set $S=PR$ and $\tilde S=P\tilde R$, where $P$ is the standard matrix permutation on $V\tp V$.
We claim that in this case the dFRT algebra associated with $\tilde S$
is a dynamical twist of $\Tg$ associated with $S$.

The FRT algebra $\Tg$ is generated by the matrix coefficients $||T^i_j||$ of
the representation $\rho$ subject to relations
\be
ST_1T_2= T_1T_2 S.
\label{FRTrel}
\ee
The $\U$-bimodule structure $\Tg$ is defined
by the formulas (\ref{bi-act}), where $x\in \U$ and $y\in \U_{op}$.

Let us consider $\Tg\tp \Lc_{op}\tp \Lc$
as an algebra over the $\Lc_{op}\tp \Lc$-bialgebroid
$\Bg_{op}\tp \Bg$, where $\Bg=\U\tp \D\Hg_\Lc$
and $\Bg_{op}$ differs from $\Bg$ by the replacement of $(\U,\Hg,\Lc)$
by $(\U_{op},\Hg_{op},\Lc_{op})$, see \cite{DM2}.
Recall from Example \ref{ex_bialg2} that $\Lc_{op}$ is an $\Hg_{op}$-base algebra.

Suppose that $\F\in \Us{\e}{\e}{2}$ is a universal dynamical
cocycle over the base $\Lc$. Then $\bar\F:=\F^{-1}\in
\U_{op}\tp \U_{op}\tp \Lc_{op}$ is a universal dynamical cocycle
over $\Lc_{op}$, \cite{DM2}. Thus $\bar\F\tp \F$ is a universal twist in the dynamical extension of the
category of $\U$-bimodules over the base algebra
$\Lc_{op}\tp \Lc$. In the standard way, it amounts to a twist of the bialgebroid
$\Bg_{op}\tp \Bg$. Denote by $\tilde \Tg'$ the corresponding twist
of the module algebra $\Tg\tp \Lc_{op}\tp \Lc$, as explained in Subsection \ref{subsecdRQG}.
\begin{propn}
The algebra $\tilde \Tg'$ is isomorphic to the dFRT algebra $\tilde \Tg$ from  Definition
\ref{dFRTa}.
\label{dFRTa'}
\end{propn}
\begin{proof}
Immediate corollary of (\ref{dyn_alg}) and (\ref{FRTrel}).
\end{proof}

Now consider the  dual  Hopf algebra $\U^*$ to $\U$ as a bimodule over $\U$.
This algebra is generated by the matrix coefficients of a fundamental representation of $\U$.
They satisfy the ordinary FRT relations plus some additional relations, e.g. the unit quantum determinant in
the $\U=\U_q(sl(n))$-case. Apply the dynamical twist $\bar \F\tp \F$ to $\U^*\rtimes( \Lc_{op}\tp \Lc)$
and construct the algebra $\tilde \U^*$ similarly to $\tilde \Tg'\simeq \tilde \Tg$ above.

Consider the matrix coefficients of the representation $\rho$ as elements of $\U^*$
arranged into the matrix $U=||U^i_j||_{i,j=1}^{\dim V}$. Observe
that $U$ is invertible in $\End(V)\tp \U^*$ and $U^{-1}=\gm_{\U^*}(U)$.
We are going to prove that the matrix $U$ is invertible when its entries are regarded as elements
of $\tilde \U^*$.
Introduce the matrix $\bar U\in \End(V)\tp \tilde \U^*$ setting
$$
\bar U:=\rho\bigl(\gm(\zt_1)\bigr) \gm_{\U^*}(U) \rho\bigl(\gm(\bar\zt_1)\bigr)\tp \bar\zt_{\Lcc}\tp\zt_{\Lcc}.
$$
Here $\gm$ is the antipode in $\U$ and $\zt$ is defined in (\ref{tz}).
Note that the component $\bar \zt_\Lcc$ lies in $\Lc_{op}$.

\begin{propn}
\label{U_inv}
The matrix $\bar U$ is the inverse to $U$ in the associative algebra $\End(V)\tp \tilde \U^*$.
\end{propn}
\begin{proof}
Let us prove that $\bar U$ is a right inverse to $U$.
If follows from definition of $*$ and $\bar U$ that
\be
\label{aux_U*U}
U*\bar U =
\bar \F_1 U \F_1 \gm(\zt_1)\gm(\F_2)\gm_{\U^*}(U)\gm \bar \F_2 \gm(\bar \zt)\tp
\bar\zt_{\Lcc}\bar\F_{\Lcc}\tp\F_\Lcc\zt_\Lcc .
\ee
Note with care that the product of $\Lc_{op}$-elements in this formula is
written in terms of $\Lc$. Also mind that we suppressed the symbol of representation $\rho$.

The right-hand side in (\ref{aux_U*U}) is equal to unit in view of the identities
\be
\left\{
\begin{array}{lll}
\F_1\gm(\zt_1)\gm(\F_2) \tp \F_\Lcc\zt_\Lcc&=&1_{\check \U}
\\
\bar \F_1\gm(\bar\F_2)\gm(\bar\zt_1) \tp \bar\zt_\Lcc\bar\F_\Lcc&=&1_{\check \U}
\end{array}
\right.
,
\ee
which hold by definition of $\zt$, cf. (\ref{tz}).

Now let us prove that $\bar U$ is a left inverse to $U$. We have
\be
\label{aux_U*U1}
\bar U*U =
\gm(\zt)\gm(\F_1)\gm_{\U^*}(U)\gm(\bar \F_1)\gm(\bar \zt)\Theta_{2'}\bar\F_2 U \F_2\Theta_2\tp
(\Theta_{1'}\tr\bar\zt_{\Lcc})\bar\F_{\Lcc}\tp\F_\Lcc(\Theta_1\tr\zt_\Lcc) .
\ee
Again, we suppressed the symbol or representation $\rho$. The product of $\Lc_{op}$-elements is written
in terms of $\Lc$.
Remark that $\Theta$-s have appeared in this expression because the first factor $\bar U$ contains elements
from $\Lc_{op}\tp \Lc$. Those elements
must be pulled through $U$ to the  right prior to applying $\F$.
This follows from associativity of $*$ and (\ref{*prod}).

The right-hand side of the equation (\ref{aux_U*U1}) is equal to unit if and only if
\be
\left\{
\begin{array}{lll}
\gm(\bar \F_1) \gm(\bar \zt_1)\Theta_2\bar \F_2\tp \bar \F_{\Lcc}(\Theta_1\tr \bar \zt_\Lcc)&=&1_{\check \U}
\\
\gm(\zt_1)\gm(\F_1)\F_2\Theta_2\tp \F_{\Lcc}(\Theta_1\tr \zt_\Lcc)&=&1_{\check \U}
\end{array}
\right.
\Leftrightarrow
\left\{
\begin{array}{lll}
\Gm(\bar\zt)&=&\vt
\\
\Gm(\zt)&=&\bar\vt
\end{array}
\right.
.
\ee
The right system of equations is satisfied by Proposition \ref{prop_vt_inv}.
\end{proof}

Let us associate with the $\Ha$-invariant matrix $\tilde S\in \End^{\tp 2}(V)\tp \Lc$
the extended dFRT algebra $\tilde \Tg_{ext}$, as in Subsection \ref{subsecEdFRTa}.
We have the following result,  as a corollary of Proposition \ref{U_inv}.
\begin{propn}
\label{T-inv}
The correspondence $T^i_j\mapsto U^i_j$, $\bar T^i_j\mapsto \bar U^i_j$ defines a bialgebroid homomorphism
from $\tilde \Tg_{ext}$ to $\tilde \U^*$.
\end{propn}
\begin{proof}
The algebra $\U^*$ is a quotient of the FRT algebra $\Tg$.
The matrix $U$ satisfies the relations (\ref{FRTrel}) in $\End(V)\tp \U^*$. Therefore
it satisfies the relations (\ref{dFRTrel}) in $\End(V)\tp \tilde\U^*$,
by Proposition \ref{dFRTa'}. Further, the matrix $\bar U$ is the inverse to $\U$, by Proposition
\ref{U_inv}.
These arguments readily imply that the relations (\ref{T-barT}) involving $\bar T$  are also  preserved
under the correspondence $T\mapsto U$, $\bar T\mapsto \bar U$.
\end{proof}
\subsection{dRE algebra via dynamical twist}
\label{subsecdREadT}
Denote by $\F_{(12)(34)}$ the image of the dynamical $(\U,\Hg,\Lc)$-twisting cocycle $\F$
under the map $\Delta\tp \Delta \colon \U\tp \U\to \mbox{\tw{\U}{\Ru}{\U}}\tp \mbox{\tw{\U}{\Ru}{\U}}$.
Since the coproduct
$\Delta\colon \U\to \mbox{\tw{\U}{\Ru}{\U}}$ is a  Hopf algebra homomorphism,
the element $\F_{(12)(34)}$ is a dynamical $(\mbox{\tw{\U}{\Ru}{\U}},\Hg,\Lc)$-twisting cocycle.
Here $\Ha$ is viewed as a subalgebra in \tw{\U}{\Ru}{\U} by the inclusion via $\Delta$.
Consider the coboundary dynamical $(\mbox{\tw{\U}{\Ru}{\U}},\Hg,\Lc)$-twist
generated by $\F\in\mbox{\tw{\U}{\Ru}{\U}}\tp \Lc$.
The composition of $\F_{(12)(34)}$ with this coboundary twist is again
a dynamical $(\mbox{\tw{\U}{\Ru}{\U}},\Hg,\Lc)$-twist whose explicit
expression is
$(\Ru^{-1}_{23}\F_{(13)(24)}\Ru_{23})\F_{(12)(34)}(\F_{12|34}\F_{34})$.
The left factor in the parentheses is obtained by the coproduct of
 $\mbox{\tw{\U}{\Ru}{\U}}$ applied to $\F$.
The dynamical  R-matrix for the triple $(\mbox{\tw{\U}{\Ru}{\U}},\Hg,\Lc)$
then reads
$$
(\F_{(34)(12)}\F_{34|12}\F_{12})^{-1}(\Ru^-_{14}\Ru_{24}\Ru^-_{13}\Ru_{23})(\F_{(12)(34)}\F_{12|34}\F_{34})=
\tilde\Ru^-_{14|2}\tilde\Ru_{24}\tilde\Ru^-_{13|24} \tilde\Ru_{23|4},
$$
where $\tilde \Ru=\bar\F_{21}\Ru\F$ and  $\tilde \Ru^-=\bar\F_{21}\Ru^-\F=\tilde\Ru_{21}^{-1}$.
Let us perform one more coboundary twist with the element $\zt_{1|2}\in \mbox{\tw{\U}{\Ru}{\U}}\tp \Lc$.
Denote by $\Upsilon$ the resulting composite dynamical twist.

Let $\Kg$ be the RE algebra associated with $S$ and generated by the matrix elements
$||K^i_j||$. It is a module algebra over \tw{\U}{\Ru}{\U} with respect to
the action
$$
(x\tp y)\btr K=\rho\bigl(\gm (x)\bigr)K\rho(y), \quad
x\tp y \in \mbox{\tw{\U}{\Ru}{\U}}.
$$
Since \tw{\U}{\Ru}{\U} coincides as an algebra with $\U\tp \U$, this action can be
expressed through right and left actions of $\U$
\be
(x\tp y)\btr A= y\rightharpoonup A \leftharpoonup \gm(x).
\label{action}
\ee
Remark that $\Kg$ coincides as a linear space with the FRT algebra $\Tg$, which is a bialgebra and admits
a Hopf pairing with $\U$. The actions $\rightharpoonup$ and $\leftharpoonup$ turns into, respectively,
left and right coregular actions of $\U$ on $\Tg$.

Consider $\K\tp \Lc$ as an algebra over the left $\Lc$-bialgebroid $(\mbox{\tw{\U}{\Ru}{\U}})\tp (\D\Ha_\Lc)$.
The bialgebroid twist induced by $\Upsilon$ transforms $\K\tp \Lc$ to a module algebra, $\tilde \Kg''$,
over the twisted bialgebroid.
Let $\tilde S$ be the dynamical matrix built by $\tilde \Ru$ as in Subsection \ref{ssDFRT-T}.
Associate with $\tilde S$ the dRE algebra $\tilde \Kg$, according to Definition \ref{defdREa}.
\begin{thm}
\label{dREbyTwist}
The algebra $\tilde \Kg''$ is isomorphic to the dRE algebra $\tilde \Kg$.
\end{thm}
\begin{proof}
The dynamical twist $\Upsilon$ converts the universal R-matrix
of \tw{\U}{\Ru}{\U} to  the dynamical R-matrix
$$
(\bar\zt_{3|412}\bar\zt_{1|2})\tilde\Ru^-_{14|2}\tilde\Ru_{24}\tilde\Ru^-_{13|24} \tilde\Ru_{23|4}(\zt_{1|234}\zt_{3|4})
\in (\mbox{\tw{\U}{\Ru}{\U}})\tp (\mbox{\tw{\U}{\Ru}{\U}})\tp \Lc.
$$
To write out the commutation relations in $\tilde\Kg''$, it is convenient
to represent this R-matrix in the form (recall that the bar means the inverse)
$$
{\cal S}'_{21}{\cal \bar S}''=(\bar\zt_{3|412}\bar\zt_{1|2})\tilde\Ru^-_{14|2}\tilde\Ru_{24}
(\zt_{3|24}\zt_{1|324})(\bar\zt_{1|324}\bar\zt_{3|24})\tilde\Ru^-_{13|24} \tilde\Ru_{23|4}(\zt_{3|4}\zt_{1|234}).
$$
where ${\cal S}'={\cal S}'_1\tp {\cal S}'_2\tp {\cal S}'_\Lc$ and
${\cal S}''={\cal S}''_1\tp {\cal S}''_2\tp {\cal S}''_\Lc$ belong to $(\mbox{\tw{\U}{\Ru}{\U}})\tp (\mbox{\tw{\U}{\Ru}{\U}})\tp \Lc$
and
\be
\begin{array}{lll}
{\cal S}':=& (\bar\zt_{1|234}\bar\zt_{3|4})\tilde\Ru_{23|4}^{-1}\tilde\Ru_{42}(\zt_{1|42}\zt_{3|142})
=
\bigl(\bar\zt_{3|4}\tilde\Ru_{23|4}^{-1}\zt_{3|24}\bigr)\bigl(\tilde\Ru_{42}\bigr)
,\\[6pt]
{\cal S}'':=&(\bar\zt_{1|234}\bar\zt_{3|4})\tilde\Ru_{23|4}^{-1}\tilde\Ru_{31|24}(\zt_{3|42}\zt_{1|324})
=
\bigl(\bar\zt_{3|4}\tilde\Ru_{23|4}^{-1}\zt_{3|24}\bigr)\bigl(\bar\zt_{3|24}\bar\zt_{1|324}\tilde\Ru_{31|24}\zt_{1|24}\zt_{3|124}\bigr).
\end{array}
\label{primes}
\ee
Then the commutation relation (\ref{dyn_alg}) for two arbitrary elements $A,B\in\Kg\subset \tilde \Kg''$ can be written as
\be
\label{relations}
({\cal S}'_1\btr B)*({\cal S}'_2\btr A)*{\cal S}'_\Lcc=
({\cal S}''_1\btr A)*({\cal S}''_2\btr B)*{\cal S}''_\Lcc,
\ee
The left equalities (\ref{r_ant1}), (\ref{r_ant2}), and (\ref{r_ant3}) can be represented as
\be
(\gm\tp \gm)(\bar\zeta_{2|1}\bar\zeta_1\tilde \Ru \zeta_{1|2}\zeta_2)
&=&
\bigl(\bar\gm(\Theta_1)\tr \tilde \Ru \bigr)\Delta(\Theta_2),
\label{0r_ant1}
\\
(\gm \tp \id)(\bar \zeta_1\tilde \Ru \zt_{1|2})
&=&
\bigl(\bar\gm(\Theta_1)\tr \tilde \Ru^{-1} \bigr)(\Theta_2\tp 1),
\label{0r_ant2}
\\
(\id\tp \gm)(\bar\zeta_2\tilde \Ru^{-1} \zt_{2|1})
&=&
\bigl(\bar\gm(\Theta_1)\tr \tilde \Ru \bigr)(1\tp \Theta_2).
\label{0r_ant3}
\ee
We substitute (\ref{primes}) into (\ref{relations}) and evaluate the action (\ref{action}).
The result will be expressed trough the left hand-side of
(\ref{0r_ant1})--(\ref{0r_ant3}). Then we replace them by the right-hand side of (\ref{0r_ant1})--(\ref{0r_ant3}).
In that way we have  for the right-hand side of (\ref{relations})
\be
\label{aux_eq1}
\bigl\{\tilde\Ru_{1'}\Theta^{(3)}_2 \rightharpoonup
A \leftharpoonup \tilde \Ru_2\gm(\Theta^{(1)}_2)\bigr\} * \tilde\Ru_{\Lcc'}  *
\bigl\{\Theta^{(4)}_2 \rightharpoonup B \leftharpoonup\tilde\Ru_1\gm(\Theta^{(2)}_2)
\tilde\Ru_{2'}\bigr\}*\bigl\{\Theta_1\tr\tilde\Ru_{\Lcc}\bigr\}.
\ee
The dynamical R-matrix has the symmetry type $(\tau,\e)$. Applying this argument to the R-matrix labelled with prime,
we rewrite (\ref{aux_eq1}) as
\be
\label{aux_eq2}
\bigl\{\Theta^{(2)}_2\tilde\Ru_{1'} \rightharpoonup
A \leftharpoonup\tilde \Ru_2\gm(\Theta^{(1)}_2)\bigr\}* \bigl\{\Theta^{(3)}_2\tr\tilde\Ru_{\Lcc'}\bigr\}
*
\hspace{2in}\nn\\
*
 \bigl\{\Theta^{(5)}_2 \rightharpoonup
B \leftharpoonup \tilde\Ru_1 \tilde \Ru_{2'}\gm(\Theta^{(4)}_2)\bigr\}* \bigl\{\Theta_1\tr  \tilde\Ru_{\Lcc}\bigr\}.
\ee
Now recall that elements of $\Kg$ and $\Lc$ permute in the algebra $\tilde \Kg$ according to the rule
$$
\la*A = \bigl(\Theta_2^{(2)} \rightharpoonup A \leftharpoonup\gm(\Theta_2^{(1)})\bigr)*(\Theta_1\tr \la)
$$
 for all $A\in \Kg$ and $\la\in \Lc$ (apply  (\ref{*prod}) and (\ref{action}) setting $x\tp y=\Delta(\Theta_2)$).
Taking this into account, pull the factor $\Theta_1\tr\tilde\Ru_{\Lcc}$ in (\ref{aux_eq2})
to the left. This will kill all the $\Theta$-s.
Thus the right-hand side of (\ref{relations}) acquires the form
\be
\tilde\Ru_{\Lcc}* ( \tilde\Ru_{1'} \rightharpoonup A \leftharpoonup\tilde\Ru_2) * \tilde\Ru_{\Lcc'} * (B \leftharpoonup\tilde\Ru_1 \tilde\Ru_{2'}).
\nn
\ee
Taking $A$ and $B$ to be the matrix coefficients $||K^i_j||$, we obtain the left-hand side of (\ref{dRE}).

For the left-hand side of (\ref{relations}) we have
\be
&&\bigl\{ \tilde\Ru_1\tilde \Ru_{2'}\rightharpoonup B \bigr\}   *
\bigl\{ \Theta^{(2)}_2\tilde\Ru_{1'}\rightharpoonup A \leftharpoonup \tilde\Ru_2\gm(\Theta^{(1)}_2)\bigr\}*
\bigl\{\Theta_1\tr\tilde\Ru_{\Lcc}\bigr\}*\tilde\Ru_{\Lcc'}=
\hspace{1.0in}\nn\\
&&\hspace{0.5in}=\bigl\{ \tilde\Ru_1\tilde \Ru_{2'} \rightharpoonup B \bigr\} * \bigl\{\tilde\Ru_{\Lcc}\bigr\} *
\bigl\{\tilde\Ru_{1'} \rightharpoonup A \leftharpoonup \tilde\Ru_2\bigr\}*\tilde\Ru_{\Lcc'}.
\nn
\ee
Taking $A$ and $B$ to be the matrix coefficients $||K^i_j||$ we obtain the right-hand side of (\ref{dRE}).

Thus we have shown that there is an algebra map $\tilde \K\to \tilde \K''$.
This map is in fact an isomorphism.  The ordinary RE algebra is the quotient
of $k\langle K^i_j\rangle$ by the RE relations. The bialgebroid twist $\Upsilon$
transforms the algebra $k\langle K^i_j\rangle\tp \Lc$, which is free over $\Lc$, into an
$\Lc$-free algebra. It is easy to see that that algebra  is isomorphic to
$\Fg(K)=k\langle K^i_j\rangle\rtimes \Lc$. Therefore $\tilde \K''$ is a quotient
of $\Fg(K)$ by the dRE relations, which are the $\Upsilon$-image of
the ordinary RE relations. Hence $\tilde \K''\simeq \tilde \K$.

\end{proof}
\subsection{A homomorphism of the dRE algebra to $\U\tp \Lc$}
The ordinary RE algebra associated with a representation of a quasitriangular
Hopf algebra $\U$ admits a homomorphism to $\U$.
Below we show that the dRE algebra has a homomorphism into $\check \U$.

First of all observe that the matrix $\tilde \Q:=\Ru_{21}\Ru\in \Us{\e}{\e}{2}$ satisfies equation
\be
\tilde\Ru_{21|3}\tilde\Q_{13}\tilde\Ru_{12|3}\tilde\Q_{23}
=
\tilde\Q_{23}\tilde\Ru_{21|3}\tilde\Q_{13}\tilde\Ru_{12|3}
\label{udRE}
\ee
in $\check \U^3$.
This equation can be derived directly from the dYBE. Another way is to interpret
the dynamical R-matrix through  the generator of a non-local representation of the braid group.
Put $\tilde Q=\rho(\tilde\Q_1)\tp \tilde\Q_2  \tp \tilde\Q_\Lcc\in \End(V)\tp \U\tp \Lc =\End(V)\tp \check\U$.
\begin{propn}
Along with the embedding $\delta\colon\Lc\to \U\tp \Lc$, the correspondence
$K^i_j \mapsto \tilde Q^i_j$ defines an algebra homomorphism
$\varphi\colon\tilde \Kg\to  \check \U$.
\end{propn}
\begin{proof}
Applying $\rho\tp \rho$ to the first two $\U$-factors in (\ref{udRE}), we
come to the equation
\be
\label{dRErep}
\tilde R_{21|3}\tilde Q_{13}\tilde R_{12|3}\tilde Q_{23}
=
\tilde Q_{23}\tilde R_{21|3} \tilde Q_{13}\tilde R_{12|3}.
\ee
This is exactly the image of the dRE equation (\ref{dRE}) under the map $\varphi$.
We must to check that the map $\varphi$ preserves the permutation
rule (\ref{Ktrans}). This is the case due to $\tilde \Q\in \Us{\e}{\e}{2}$.
\end{proof}

\begin{remark}
Equation (\ref{dRErep}) is equivalent to the dRE considered in \cite{FHS} for
the case of abelian base.
Indeed, multiply equation (\ref{udRE}) by $\zt_{1|23}\zt_{2|3}\zt_{3}$ from the right and by
$\bar\zt_{1|23}\bar\zt_{2|3}\bar\zt_{3}$ from the left, then apply $\Gm_{123}$.
One can check, using Proposition \ref{Rmatsym} (compare this also
with  Example \ref{autdYBE}), that this operation converts (\ref{udRE}) into
\be
\tilde\Q_{23|1}\tilde\Ru_{21}\tilde\Q_{13|2}\tilde\Ru_{12}
=
\tilde\Ru_{21}\tilde\Q_{13|2}\tilde\Ru_{12}\tilde\Q_{23|1}.
\label{udRE1}
\ee
Passing to a representation we come to the dRE that appeared in \cite{FHS}.
\end{remark}

\section{Dynamical trace and center of the dRE algebra}
\label{secDTC}
\subsection{Dynamical trace}
In this section we introduce dynamical trace (in the twisted setting) and use
it to construct central elements of the dRE algebra.

Set $\mu:=\Ru_1\gm(\Ru_2)\in \U$. The conjugation with $\mu$ implements the inverse squared antipode in $\U$,
\cite{Dr2}.
Let $(V,\rho)$ be an $\U$-module and let $\A$ be a module over the left bialgebroid $\Lc\rtimes\Ha$.
Then $\A$ is an $\Lc$-bimodule. The right $\Lc$-action is expressed through the left $\Lc$-action
and the action $\btr$ of $\Ha$ by $\la a = (\la^{(1)}\btr a) \la^{[\infty]}$, where $a\in \A$ and $\la\in \Lc$.
We will say that a matrix $X\in \End(V)\tp \A$ is $\Ha$-equivariant if its coefficients
form the dual module to $\End(V)$, namely $h\btr X =\rho\bigl(\gm(h^{(1)})\bigr)X\rho(h^{(2)})$, for all $h\in \Ha$.

To simplify further formulas we will suppress the symbol of representation $\rho$ whenever possible.
\begin{definition}
Let $X$ an $\Ha$-invariant matrix with coefficients in $\A$.
Dynamical trace $\Tr_d(X)$ is the element of $\A$ defined  by
\be
\Tr_d(X):=\Tr(\mu\vt X \zt).
\label{d_trace}
\ee
\end{definition}
Formula (\ref{d_trace}) gives a generalization to the dynamical case of the
 quantum trace for a quasitriangular Hopf algebra and its twist.
Let us prove the following important properties of this generalization.
\begin{lemma}
\label{prop_tr}
The dynamical trace satisfies the following equalities:
\be
\Tr\bigl(X \zt\bar \Gm(\bar\zt)\mu \bigr)=\Tr_d(X)=\Tr\bigl(\mu \Gm(\bar \vt)\vt X\bigr).
\label{d_trace1}
\ee
\end{lemma}
\begin{proof}
We will prove only the right equality in (\ref{d_trace1}); the left one is verified similarly.
If the coefficients of the matrix $\rho(\zt)$ were scalars, we could have moved $\zt$ in  $\Tr(\mu\vt X \zt)$
freely to the left-most position. However we should take into account how elements from $\Lc$ permute with
the coefficients of the matrix $\rho(\vt)X$. The latter is $\Ha$-invariant, and
we can transform $\Tr(\mu\vt X \zt)$ to
\be
(\Tr\tp \id)\Bigl(\bar\gm(\Theta_{2'})\zt_1\mu\vt_1\Theta_2 \tp \vt_\Lcc\bigl((\Theta_1\Theta_{1'})\tr\zt_\Lcc)\bigr)
\sum_{i,j}e^i_j\tp X^j_i \Bigr).
\ee
As an element from $\Us{\e}{\e}{}$, $\vt$ commutes with $\delta(\Lc)$. Therefore
the factor before the sum is equal to
$
\bar\gm(\Theta_{2'})\zt_1\mu\Theta_2\vt_1 \tp \bigl((\Theta_1\Theta_{1'})\tr\zt_\Lcc)\bigr)\vt_\Lcc \in \check \U.
$
Pulling $\mu$ to the left we bring this expression to $\mu \Gm^2(\zt)\vt=\mu \Gm(\bar \vt)\vt$.
This proves the right equality in (\ref{d_trace1}).
\end{proof}
\begin{lemma}
\label{lemma_trace_aux}
Let $X\in \End(V)\tp \A$ be an $\Ha$-invariant matrix.
Then
\be
\Tr_1(\mu_1\vt_{1|2} \tilde \Ru^{-1}_{12} X_1\tilde \Ru_{12}\zt_{1|2})
=
1_{\End(V)}\tp \Tr_d(X)
=
\Tr_1(\mu_1\vt_{1|2} \tilde \Ru_{21} X_1\tilde \Ru^{-1}_{21}\zt_{1|2}).
\label{aux_trace}
\ee
\end{lemma}
\begin{proof}
Using (\ref{r_ant2}) and (\ref{r_ant3}) we rewrite (\ref{aux_trace}) as
\be
\Tr_1\Bigr(\mu_1\Gm_1(\tilde \Ru_{12})\vt_{1}  X_1\zt_{1}\bar\Gm_1(\tilde \Ru^{-1}_{12})\Bigl)
=
1\tp \Tr_d(X)
=
\Tr_1\Bigr(\mu_1\Gm_1(\tilde \Ru^{-1}_{21})\vt_1  X_1\zt_{1}\bar\Gm(\tilde \Ru_{21})\Bigl).
\label{aux_trace1}
\ee
We will check only the left equality in (\ref{aux_trace1}).
We will see that it holds for an arbitrary matrix $\tilde R$, therefore
our proof will be valid upon replacement
$\tilde R\leftrightarrow \tilde R^{-1}_{21}$.

Pull the matrix $\bar\Gm_1(\tilde \Ru^{-1}_{12})$ to the left under the trace
taking into account how elements from $\Lc$ permutes with coefficients of the matrix
$\rho(\vt)  X\rho(\zt)$ (cf. the proof of Proposition \ref{prop_tr}).
The latter is $\Ha$-invariant, therefore
the expression under $\Tr_1$ can be brought to the form
$$
\Bigl(\bar\gm(\Theta_{2''})\Theta_{2'}\bar\gm(\tilde \Ru^{-1}_{1'})\mu\gm(\bar\Theta_{2}\tilde \Ru_1)\Theta_{2'''}
\tp
\tilde \Ru_2\tilde \Ru^{-1}_{2'}\tp (\bar\Theta_{1}\tr\Ru_{\Lcc})
\bigl((\Theta_{1'''}\Theta_{1''}\Theta_{1'})\tr\Ru^{-1}_{\Lcc'}\bigr)\Bigr)(\vt_{1}  X_1\zt_{1})
$$
(mind that we suppress $\rho$).
We will prove the lemma if show that the left factor within the big parentheses equals $\mu_1\in \check \U^2$.
First of all notice that the term  $\bar\gm(\Theta_{2''})\Theta_{2'}\tp \Theta_{1''}\Theta_{1'}$
cancels. Then pull $\mu$ to the right and employ the fact that conjugation with $\mu$ implements $\gm^{-2}$.
The resulting expression in the big brackets will be
$$
\bar\gm(\tilde \Ru^{-1}_{1'})\bar\gm(\bar\Theta_{2}\tilde \Ru_1)\bar\gm(\bar\Theta_{2'''})\mu
\tp
\tilde \Ru_2\tilde \Ru^{-1}_{2'}\tp (\bar\Theta_{1}\tr\tilde\Ru_{\Lcc})(\bar\Theta_{1'''}\tr\tilde\Ru^{-1}_{\Lcc'}),
$$
which is equal to
$$
\bar\gm(\tilde \Ru_1\tilde \Ru^{-1}_{1'})\bar\gm(\bar\Theta_{2})\mu
\tp
\tilde \Ru_2\tilde \Ru^{-1}_{2'}\tp \bar\Theta_{1}\tr(\tilde\Ru_{\Lcc}\tilde\Ru^{-1}_{\Lcc'})=\mu_1 .
$$
Thus the leftmost  term in (\ref{aux_trace}) equals $\Tr_1(\mu_1\vt_{1}  X_1\zt_{1})$, which proves the left quality.
\end{proof}
\begin{remark}
The dynamical trace can be defined for any (dualizable) solution to dYBE, without the dynamical twist assumption.
Indeed, such solution (multiplied by permutation) is just a Yang-Baxter operator in the dynamical category,
which is rigid, cf. Remark \ref{rigid}.
For the theory of dualizable Yang-Baxter operators in rigid categories, consult e.g. \cite{JS}.
\end{remark}
\subsection{On the center of the dRE algebra}
We will apply the dynamical trace to constructing central elements of the algebra $\tilde \Kg$.
We set $\A=\tilde \Kg$ and notice that any power $K^n$, $n\in \N$, of the matrix $K$
is $\Ha$-invariant.
\begin{thm}
\label{center_trace}
Let $K$ be the matrix of the generators of the dRE algebra $\tilde \Kg$.
Then for all  $n\in \N$ the elements $\Tr_d (K^n)$ belong to the center of $\tilde \Kg$.
\end{thm}
\begin{proof}
The dRE implies the following equation for any matrix power $K^n$:
$$
\tilde R^{-1}_{12}K^n_1\tilde R_{12}K_2=K_2\tilde R_{21}K^n_1\tilde R^{-1}_{21}.
$$
Multiply this equation by $\mu_1\vt_{1|2}$ from the left and by $\zt_{1|2}$ from
the right. Since $K$ is $\Ha$-invariant,
the matrices $\K_2$ and $\xi_{1|2}$ commute for any $\xi\in \End(V)\tp \Lc$.
Now apply Lemma \ref{lemma_trace_aux}.
\end{proof}

Let us look at how the center $\Kg$ of the RE algebra $\Kg$ transforms under the dynamical twist.
It is known that all $\U$-invariant elements  of $\Kg$  lie in the center.
Let us say a few words about the nature of this phenomenon.
The $\U$-module structure on $\Kg$ is induced by the Hopf algebra
embedding $\U\stackrel{\Delta}{\longrightarrow}\mbox{\tw{\U}{\Ru}{\U}}$.
The algebra $\Kg$ is commutative in the category of \tw{\U}{\Ru}{\U}-modules.
One tensor factor of the universal R-matrix of \tw{\U}{\Ru}{\U} lies in
the image $\Delta(\U) \subset\mbox{\tw{\U}{\Ru}{\U}}$.
Therefore the universal R-matrix reduces to unit whenever that  factor
acts on a $\U$-invariant. Thus any such element commutes with entire $\Kg$.

Similar effect takes place when a quasitriangular Hopf algebra
is replaced by a quasitriangular bialgebroid.
Let us specialize our consideration to the dRE algebra.
Identify $\Kg$ with a subspace $\Kg\tp 1\subset\tilde \Kg$ and
consider $\tilde \Kg$ as a natural module $\simeq\Kg\tp \Lc$ over the associative algebra $(\U\tp\U)\tp \Lc$.
\begin{propn}
\label{center_invariance}
For any $\U$-invariant element $a\in \Kg$ the element
$$\Bigl(\bigl(\zt\bar\Gm(\bar\zt)\bigr)_1\rightharpoonup a\Bigr)\tp\bigl(\zt\bar\Gm(\bar\zt)\bigr)_\Lcc$$
is central in $\tilde \Kg$.
\end{propn}
\begin{proof}
Consider the algebra $\tilde \Kg'$ obtained from $\Kg$ by the dynamical twist $\F_{(12)(34)}$, as
in Subsection \ref{subsecdREadT}.
Multiplication in $\tilde \Kg'$ (as well as in $\tilde \Kg$) makes it a free right $\Lc$-module
$\simeq\Kg\Lc$. Therefore, to prove centrality of an element, it suffices to show that it commutes
with elements from $\Kg$ and $\Lc$.

Every $\U$-invariant element  $a$  from $\Kg$ is $\Ha$-invariant  and therefore commutes with $\Lc$
in the algebra  $\tilde \Kg'$, cf. (\ref{*prod}). Let us show
 that it also commutes with elements from $\Kg\subset \tilde \Kg'$.
 Indeed,
the permutation of $a$ with elements from $\Kg$ is controlled by the dynamical R-matrix,
cf. (\ref{dyn_alg}). In our case that dynamical R-matrix
is obtained from the R-matrix of \mbox{\tw{\U}{\Ru}{\U}} by dynamical twist
$\F_{(13)(24)}$. Both $\Ru$ and $\F_{(13)(24)}$ cancel on $a$. Therefore
$a$ commutes with $\Kg$ and is thus central in $\tilde \Kg'$.

By Theorem \ref{dREbyTwist}, the algebras $\tilde \Kg$ and $\tilde \Kg'$  are module algebras over bialgebroids that
differ from each other by the coboundary twist
$
\F_{12}\zt_{1|2}\in \mbox{\tw{\U}{\Ru}{\U}}\tp \Lc \subset \mbox{\tw{\U}{\Ru}{\U}}\tp \D\Ha_\Lc
$.
Therefore $\tilde \Kg'$ and $\tilde \Kg$ are isomorphic as associative algebras,
and that isomorphism is implemented via the element $\F_{12}\zt_{1|2}$:
$$
\tilde \Kg'\ni b\mapsto (\F_{12}\zt_{1|2})^{-1}\btr b\in \tilde \Kg.
$$
Here we denote the action of $\mbox{\tw{\U}{\Ru}{\U}}\tp \D\Ha_\Lc$ by the same symbol as the
action of \tw{\U}{\Ru}{\U}.
To finish the proof, one should evaluate this map taking into account
$$(x\tp y)\btr a= \bigl(y\gm(x)\bigr)\rightharpoonup a$$
for any $\U$-invariant element $a\in \Kg\subset \tilde \Kg$ and any $x\tp y\in \mbox{\tw{\U}{\Ru}{\U}}$.
\end{proof}
\begin{remark}
Proposition \ref{center_invariance} tells us that the center of $\tilde \Kg$ is not less than
the center of $\Kg$.
However, Proposition \ref{center_invariance} does not directly imply Theorem \ref{center_trace}.
\end{remark}

\vspace{6pt}
{\bf Acknowledgement.}
A. Mudrov is grateful to the Max-Planck Institut f\"{u}r Mathematik for hospitality.
He appreciates numerous useful and interesting discussions with  P. Pyatov.
He also thanks D. Gurevich for drawing attention to the problem of dynamical trace and
center of the dynamical reflection equation algebra.

\end{document}